\documentclass[10pt]{amsart}
\usepackage{epsfig}
\usepackage{amsmath, amssymb, euscript,  enumerate}
\usepackage{bbm}
\usepackage{color}
\usepackage{verbatim}
\newcommand{\R}{{\mathbb R}}

\newcommand{\sgn}{{\text{\rm sign}}}

\newcommand{\diag}{{\text{\rm diag}}}

\newcommand{\ap}{\alpha}             

\newcommand{\gm}{\gamma}             \newcommand{\Gm}{\Gamma}
\newcommand{\dt}{\delta}             
\newcommand{\vep}{\varepsilon}
\newcommand{\zt}{\zeta}
\newcommand{\Th}{\Theta}

\newcommand{\kp}{\kappa}

\newcommand{\ld}{\lambda}            
\newcommand{\sm}{\sigma}             
\newcommand{\vp}{\varphi}
\newcommand{\om}{\omega}             
            \newcommand{\iy}{\infty}
            
\newcommand{\f}{\frac}             \newcommand{\el}{\ell}

\newcommand{\fA}{{\mathfrak A}}

\newcommand{\fG}{{\mathfrak G}}

\newcommand{\fL}{{\mathfrak L}}

\newcommand{\fS}{{\mathfrak S}}

\newcommand{\fU}{{\mathfrak U}}

\newcommand{\fe}{{\mathfrak e}}

\newcommand{\fr}{{\mathfrak r}}

\newcommand{\BC}{{\mathbb C}}

\newcommand{\BH}{{\mathbb H}}

\newcommand{\BN}{{\mathbb N}}

\newcommand{\BQ}{{\mathbb Q}}
\newcommand{\BR}{{\mathbb R}}

\newcommand{\BT}{{\mathbb T}}

\newcommand{\BZ}{{\mathbb Z}}

\newcommand{\cA}{{\mathcal A}}

\newcommand{\cF}{{\mathcal F}}
\newcommand{\cG}{{\mathcal G}}

\newcommand{\cM}{{\mathcal M}}

\newcommand{\cO}{{\mathcal O}}

\newcommand{\cS}{{\mathcal S}}

\newcommand{\cU}{{\mathcal U}}

\newcommand{\la}{\langle}          \newcommand{\ra}{\rangle}
\newcommand{\s}{\setminus}         
            \newcommand{\e}{\eta}

    \newcommand{\ds}{\displaystyle}

\newcommand{\bbx }{\bold x}          \newcommand{\bby }{\bold y}
\newcommand{\bbz }{\bold z}          \newcommand{\bba }{\bold a}
\newcommand{\bbb }{\bold b}          \newcommand{\bbc }{\bold c}
\newcommand{\bbd }{\bold d}
          \newcommand{\bbu }{\bold u}
\newcommand{\bbs }{\bold s}
\newcommand{\bbv }{\bold v}          
\newcommand{\bbh }{\bold h}         \newcommand{\bbn }{\bold n}
\newcommand{\bap }{\boldsymbol\alpha}
\newcommand{\bphi }{\boldsymbol\phi} 
\newcommand{\bbm }{\bold m}         \newcommand{\bbw }{\bold w}
\newcommand{\bxi }{\boldsymbol\xi}   \newcommand{\bta }{\boldsymbol\eta}
\newcommand{\bG }{\bold G}
\newcommand{\bzt }{\boldsymbol\zeta} \newcommand{\pf }{\noindent{\it Proof. }}

\newcommand{\rB }{{\text{\rm B}}}   \newcommand{\rC }{{\text{\rm C}}}
\newcommand{\rA }{{\text{\rm A}}}   \newcommand{\rD }{{\text{\rm D}}}
\newcommand{\rL }{{\text{\rm L}}}   \newcommand{\rSp }{{\text{\rm Sp}}}
\newcommand{\rSL }{{\text{\rm SL}}}
   \newcommand{\rR}{{\text{\rm R}}}

\newcommand{\rF}{{\text{\rm F}}}  \newcommand{\rS}{{\text{\rm S}}}
\newcommand{\rI}{{\text{\rm I}}}  \newcommand{\rT}{{\text{\rm T}}}
\newcommand{\rU}{{\text{\rm U}}} 
  \newcommand{\rM}{{\text{\rm M}}}
\newcommand{\rO}{{\text{\rm O}}} \newcommand{\rW}{{\text{\rm W}}}
 
\newcommand{\RE}{{\text{\rm Re}}} 
\newcommand{\lcm}{{\text{\rm lcm}}}
\newcommand{\Ima}{{\text{\rm Im}}}

\newtheorem{thm}[subsection]{Theorem}
\newtheorem{lemma}[subsection]{Lemma}
\newtheorem{cor}[subsection]{Corollary}
\newtheorem{remark}[subsection]{Remark}
\newtheorem{prop}[subsection]{Proposition}

\newtheorem{definition}[subsection]{Definition}

\newtheorem{pro}{Property}
\numberwithin{equation}{section}

\begin{document}

\title[Distribution of integral lattice points in an ellipsoid]
{Distribution of integral lattice points in an ellipsoid with a
diophantine center}

\author{Jiyoung Han}
\address{Department of Mathematics, Seoul National University, Seoul 151-
747, Korea} \email{madwink0@snu.ac.kr}
\author{Hyunsuk Kang}
\address{Department of Mathematics, Korea Institute for Advanced Study, Seoul 130-722, Korea}
\email{h.kang@kias.re.kr}
\author{Yong-Cheol Kim}
\address{Department of Mathematics Education, Korea University, Seoul
136-701, Korea} \email{ychkim@korea.ac.kr}
\author{Seonhee Lim}
\address{Department of Mathematics, Seoul National University, Seoul 151-
747, Korea} \email{slim@snu.ac.kr}

\maketitle

\begin{abstract}
We evaluate the mean square limit of exponential sums related with a
rational ellipsoid, extending a work of Marklof. Moreover, as a
result of it, we study the asymptotic values of the normalized
deviations of the number of lattice points inside a rational
ellipsoid and inside a rational thin ellipsoidal shell.
\end{abstract}

\section{Introduction}

Let $\fS^+_n(\BZ)$ be the family of all positive definite and
symmetric $n\times n$ matrices with integral components. For
$\rM\in\fS^+_n(\BZ)$, we consider the quadratic form $Q_{\rM}$
defined by $Q_{\rM}(\bbx)=\la\rM\bbx,\bbx\ra$ for $\bbx\in\BR^n$ and
the corresponding ellipsoid
$$E^{\rM}_R=\{\bbx\in\BR^n:Q_{\rM}(\bbx)\le R^2\}.$$

Our interests are focused on the distribution of integral lattice
points inside $E^{\rM}_R$ as the radius $R$ tends to infinity. In
place of the ellipsoid centered at the origin, we consider the
ellipsoid with a diophantine center of type $\kappa$ as defined
below.

\begin{definition} A vector $\bap\in\BR^n$ is said to be {\rm of
diophatine type $\kappa$}, if there exists a constant $C>0$ such
that $$\bigl|\bap-\f{\bbm}{q}\bigr|>\f{C}{q^{\kappa}}$$ for all
$\bbm\in\BZ^n$ and $q\in\BN$.
\end{definition}

The smallest possible value of $\kappa$ is $1+1/k$ in the above
definition. In this case, $\bap$ is called {\it badly approximable}
(see \cite{S}).

We now consider the counting function $N_{\rM}$ in the ellipsoid
with a diophantine center of type $\kappa$ introduced in \cite{b-b}.
Let $\mathbbm{1}_{B}$ be the characteristic function of the unit
open ball $B=B_1$ in $\BR^n$ and $\mathbbm{1}_{E^{\rM}}$ be the
characteristic function of the ellipsoid $E^{\rM}=E^{\rM}_1$
corresponding to $\rM\in\fS^+_n(\BZ)$. We write the number
$N_{\rM}(t):=\sharp\bigl(\BZ^n\cap E^{\rM}_t(\bap)\bigr)$ of lattice
points inside the ellipsoid
$E^{\rM}_t(\bap)=\{\bbx\in\BR^n:Q_{\rM}(\bbx-\bap)\le t^2\}$ of
radius $t>0$ centered at a diophantine vector $\bap\in\BR^n$; that
is to say,
$$N_{\rM}(t)=\sum_{\bbm\in\BZ^n}\mathbbm{1}_{E^{\rM}}\biggl(\f{\bbm-\bap}{t}\biggr).$$

In this paper, we investigate the asymptotics of the following
deviations: the normalized deviation $F_{\rM}(t)$ of $N_{\rM}(t)$
defined by
\begin{equation}F_{\rM}(t):=\f{N_{\rM}(t)-|E^{\rM}|\,t^n}{t^{(n-1)/2}}
\,\,\,\text{ as $t\to\iy$ }\end{equation} and the normalized
deviation $S_{\rM}(t,\e)$ of the number of lattice points inside the
spherical shell between the elliptic spheres of radii $t+\e$ and $t$
given by
\begin{equation}S_{\rM}(t,\e):=\f{N_{\rM}(t+\e)-N_{\rM}(t)-|E^{\rM}|((t+\e)^n-t^n)}{\sqrt\e\,t^{(n-1)/2}}
\end{equation} as $t\to\iy$ and $\e\to 0$, where
$|E^{\rM}|$ denotes the volume of the ellipsoid $E^{\rM}$.

Let $\mu\in C^{\iy}_c(\BR)$ be a nonnegative function which is
supported on $(0,\iy)$ and $\int_{\BR}\mu(t)\,dt=1$. Then we
consider an averaging operator $\la\,\cdot\,\ra_T$ on $L^1_{\rm
loc}(\BR)$ defined by the smooth measure induced by $\mu$ as
follows; for $f\in L^1_{\rm loc}(\BR)$,
$$\la f\ra_T:=\la f\ra_{\mu,T}:=\int_{\BR}f(t)\mu_T(t)\,dt\,\,\,\text{ for
$\mu_T(t)=\f{1}{T}\,\mu\biggl(\f{t}{T}\biggr)$ and $T>0$.}$$

Throughout the paper, we use the following notation. For two
quantities $a$ and $b$, we write $a\ll b$ (resp. $a\gg b$) if there
is a constant $C>0$ ( independent of $t\ge 1$, $T\ge 1$, $\e\ge 0$
and the smoothing parameter $K\ge 1$, but possibly depending on the
dimension $n$ and the matrix $\rM\in\fS^+_n(\BZ)$ ) such that $a\le
C b$ (resp. $b\le C a$).

The asymptotic result for the expectation values is relatively easy,
comparing with that of the variances.

\begin{thm} For any $\bap\in\BR^n$ with $n\ge 2$ and $\rM\in\fS^+_n(\BZ)$, we have that

\noindent $(a)$ $\ds\lim_{T\to\iy}\la F_{\rM}\ra_T=0$, and $(b)$
$\ds\lim_{T\to\iy}\la S_{\rM}(\,\cdot\,,\e)\ra_T=0$ provided that
$T^{-L}<\e\ll 1$ for some $L>0$.
\end{thm}

As we shall see later, this shift of the center from the origin to a
diophantine vector makes exponential sums appear in the expansion of
the counting function. For $p\in\BN$, $\rM\in\fS^+_n(\BZ)$ and
$\bap\in\BR^n$, we denote by
\begin{equation}\fr[\rM,\bap](p)=\sum_{\substack{\bbm\in\BZ^n \\ Q_{\rM}(\bbm)=p}}
\fe(\bbm\cdot\bap)
\end{equation} where $\fe(t)=e^{-2\pi i t}$ for $t>0$. A crucial ingredient for our problem is the asymptotics
of the sum
\begin{equation}\rR[\rM,\bap](N)=\sum_{p=1}^N|\fr[\rM,\bap](p)|^2\,\,\text{
as $N\to\infty$, }
\end{equation} which was studied in \cite{M3} under the condition that $\bap$ is diophantine of type
$\kappa <\f{n-1}{n-2}$ and $\rM=\rI_n$.

In case that $\rM=\rD_{\bba}:=\diag(a_1,\cdots,a_n)$ is a diagonal
matrix with diagonal entries $a_1,\cdots,a_n$ for any
$\bba=(a_1,\cdots,a_n)\in\BN^n$, we evaluate the mean square limit
of exponential sums corresponding to the ellipsoid $E^{\rD_{\bba}}$.

\begin{thm} Let $\bap\in\BR^n$ be
a vector of diophantine type $\kappa<\f{n-1}{n-2}$ such that the
components of the vector $(\bap,1)\in\BR^{n+1}$ are linearly
independent over $\BQ$. If $\,\bba\in\BN^n$, then we have that
\begin{equation}\lim_{N\to\iy}\f{1}{N^{n/2}}\,\rR[\rD_{\bba},\bap](N)
=|E^{\rD_{\bba}}|.\end{equation}
\end{thm}

The above result works for $n=1$ without the diophantine condition;
as a matter of fact, by the equidistribution of the sequence $m\bap$
modulo $1$, it turns out that the limit is $2$ for any irrational
$\bap\in\BR$. In case that $n\ge 2$ and $\rD_{\bba}=\rI_n$, the
diophantine condition is crucial, since the mean square value might
diverge for any rational $\bap\in\BR^n$ (see \cite{b-b}). Theorem
1.4 was proved by Bleher and Dyson \cite{b-d1}, in the case that
$n=2$ and $\rD_{\bba}=\rI_2$. Theorem 1.3 was shown by Marklof
\cite{M3} when $\rM=\rI_n$. Theorem 1.3 is an extension of Marklof's
work \cite{M3} to arbitrary rational ellipsoids with no rotation
(refer to Section 6 for the detailed argument). His work heavily
relies on representation theory and ergodic theory (see \cite{LV}).
Our proof follows the approach taken by Marklof in \cite{M1,M2,M3}.

One may generalize Theorem 1.3 to the case $\rM\in\fS^+_n(\BZ)$
under the same condition on $\ap$, though the validity of this
extension is not claimed in this work.

\begin{pro} Let $\bap\in\BR^n$ be a vector of
diophantine type $\kappa<\f{n-1}{n-2}$ such that the components of
the vector $(\bap,1)$ are linearly independent over $\BQ$. Then we
say that {\rm Property 1 holds} if one has
\begin{equation}\lim_{N\to\iy}\f{1}{N^{n/2}}\,\rR[\rM,\bap](N)
=|E^{\rM}|\,\,\,\text{ for $\rM\in\fS^+_n(\BZ)$ }.\end{equation}
\end{pro}

Should this property hold, we obtain the results for the variances
as stated below. We note that Theorem 1.3 is a special case in which
Property 1 holds.

\begin{thm} Let $\bap\in\BR^n$ be a
vector of diophantine type $\kappa<\f{n-1}{n-2}$ $($$n\ge 2$$)$ such
that the components of the vector $(\bap,1)\in\BR^{n+1}$ are
linearly independent over $\BQ$. For $\rM=\rD_{\bba}$ with
$\bba\in\BN^n$, the series
\begin{equation}\fA[\widehat\rM,\bap]:=(\det\rM)^{\f{n-1}{2}}\sum_{p=1}^{\iy}|\fr[\widehat\rM,\bap](p)|^2
p^{-\f{n+1}{2}}\end{equation} converges where $\widehat\rM$ is the
adjugate matrix of $\rM$, and
\begin{equation}
\lim_{T\to\iy}\la|F_{\rM}|^2\ra_T=\f{1}{2\pi^2}\,\fA[\widehat\rM,\bap].\end{equation}
Moreover, the same results follow for $\rM\in\fS^+_n(\BZ)$, if
{\,\rm Property 1} holds.
\end{thm}

\begin{thm} Let $\bap\in\BR^n$ be a
vector of diophantine type $\kappa<\f{n-1}{n-2}$ $($$n\ge 2$$)$ such
that the components of the vector $(\bap,1)\in\BR^{n+1}$ are
linearly independent over $\BQ$. For $\rM=\rD_{\bba}$ with
$\bba\in\BN^n$ and $\e\gg T^{-\gm}$ for some $\gm\in(0,1)$, we have
that
\begin{equation}
\lim_{T\to\iy}\la|S_{\rM}(\,\cdot\,,\e)|^2\ra_T=n|E^{\rM}|\,\,\,\text{
as $\e\to 0$. }
\end{equation} Moreover, the same result follows for $\rM\in\fS^+_n(\BZ)$, if
{\,\rm Property 1} holds.
\end{thm}

In case that $n\ge 2$ and $\rM=\rI_n$, Kang and Sobolev \cite{KS}
proved Theorem 1.2 and the above two theorems by applying the result
\cite{M3} of Marklof. Our results on those two theorems are
generalizations of the standard ball case to some ellipsoid case
(rotation being possibly allowed under Property 1). These results
are established by employing Theorem 1.3 or Property 1.

The paper is organized as follows. In Section 2, preliminaries on
the counting functions using Bessel functions are given. In Section
3, Theorem 1.3, the ergodic part of the paper is explained extending
the work of Marklof on Euclidean balls to ellipsoids with positive
integer coefficients. With Property 1 assumed to be hold, we give
the proofs of Theorem 1.2 and Theorem 1.4 in Section 4, and that of
Theorem 1.5 in Section 5.

\noindent{\bf Remark.} If Property 1 holds, then our results on
$\rM\in\fS^+_n(\BZ)$ to be obtained in this paper can be extended to
the case $\rM\in\fS^+_n(\BQ)$. Indeed, any $\rM\in\fS^+_n(\BQ)$ can
be written as a $n\times n$ matrix with rational components
$p_{ij}/q_{ij}$, where $p_{ij}\in\BZ$ and $q_{ij}\in\BN$ are
relatively prime. Let $c\in\BN$ be the least common multiple of all
$q_{ij}$'s. Then we see that $\rM=\f{1}{c}\,\rM_c$ for some
$\rM_c\in\fS^+_n(\BZ)$ and $E^{\rM}_R=E^{\rM_c}_{\sqrt{c}R}$ for
$R>0$. Hence we can extend all the results we got in the above to
the case $\rM\in\fS^+_n(\BQ)$, provided that Property 1 holds. Since
Theorem 1.3 is a special case of Property 1, this implies that our
results can be extended to the case that $\rM=\rD_{\bba}$ for
$\bba\in\BQ_+^n$.

If $\rM\in\fS^+_n(\BR)$ has only rational eigenvalues, then Property
1 is a rotated version of Theorem 1.3. Hence it would be interesting
to classify all $\rM\in\fS^+_n(\BR)$ for which Property 1 holds.

\section{Preliminaries}

For any $\rM\in\fS^+_n(\BZ)$, we consider the quadratic form
$Q_{\rM}$ defined by $Q_{\rM}(\bbx)=\la\rM\bbx,\bbx\ra$. Then there
is a constant $C>0$ such that $Q_{\rM}(\bbx)\ge C\,|\bbx|^2$ for any
$\bbx\in\BR^n$. By the diagonalization process, there is an
orthogonal matrix $\rO\in\cO_n(\BR)$ such that $\rM=\rO\rD\rO^{-1}$
where $\rD$ is the diagonal matrix with diagonal entries given by
positive eigenvalues $\ld_k>0$ of $\rM$. We define the square root
matrix $\sqrt\rM$ of $\rM$ by $\sqrt\rM=\rO\sqrt\rD\rO^{-1}$ where
$\sqrt\rD=\diag(\sqrt{\ld_1},\cdots,\sqrt{\ld_n})$, and define a
norm $|\cdot|_M$ on $\BR^n$ by $|\bbx|_M=\sqrt{\la\rM\bbx,\bbx\ra}$
for any $\bbx\in\BR^n$. Then we have that
$|\bbx|_M=|\sqrt\rM\,\bbx|$ for any $\bbx\in\BR^n$. Let $\fL$ be the
lattice such that $\fL=\rT(\BZ^n)$ where $\rT$ is the linear
transformation on $\BR^n$ given by the matrix $\sqrt\rM$. Then we
easily see that its dual lattice is $\fL^*=\rT^*(\BZ^n)$ where
$\rT^*$ is the linear transformation on $\BR^n$ represented by the
inverse matrix $\sqrt\rM^{\,-1}$ of $\sqrt\rM$.

Next we obtain an elementary lemma which facilitates our estimates.

\begin{lemma} Let $\rM\in\fS^+_n(\BZ)$ be given. Then we have that

$(a)$ $\widehat\rM\in\fS^+_n(\BZ)$, where $\widehat\rM$ is the
adjugate matrix of $\rM$,

$(b)$ $|\rT(\bbx)|=|\bbx|_{\rM}$ and
$|\rT^*(\bbx)|=(\det\rM)^{-1/2}|\bbx|_{\widehat{\rM}}$ for any
$\bbx\in\BR^n$.
\end{lemma}

\pf (a) By the definition of $\widehat\rM$, we see that it is
symmetric and its components are integral values. Since
$\det(\rM)^{-1}\widehat\rM=\rM^{-1}=\rO\rD^{-1}\rO^{-1}$ is positive
definite, so is $\widehat\rM$.

(b) The first part is trivial. Since
$\rM^{-1}=(\det\rM)^{-1}\widehat\rM$, we have that
\begin{equation*}\begin{split}|\rT^*(\bbx)|&=|\sqrt\rM^{\,-1}\bbx|=\sqrt{\la
\rM^{\,-1}\bbx,\bbx\ra}=\f{1}{(\det\rM)^{1/2}}\sqrt{\la
\widehat{\rM}\bbx,\bbx\ra}=\f{|\bbx|_{\widehat\rM}}{(\det\rM)^{1/2}}
\end{split}\end{equation*}
for any $\bbx\in\BR^n$. Hence we complete the proof. \qed

Let $\cS(\BR^n)$ be the Schwartz space. For $f\in\cS(\BR^n)$, its
Fourier transform $\widehat f$ on $\BR^n$ is defined by
$$\widehat f(\bxi)=\int_{\BR^n}e^{-2\pi i \bbx\cdot\bxi}f(\bbx)\,d\bbx.$$
To define a regularized version of $N_{\rM}(t)$ for
$\rM\in\fS^+_n(\BZ)$, we use a positive cut-off radial function
$\vp\in\cS(\BR^n)$ ( i.e. $\vp(\bbx)=\vp(|\bbx|)$ ) so that $\vp$ is
supported in some ball $B_R=\{\bbx\in\BR^n:|\bbx|\le R\}$ and
\begin{equation}\widehat\vp(0)=\int_{\BR^n}\vp(\bbx)\,d\bbx=1.\end{equation} For any $\vep>0$,
let $\Phi_{\vep}(\bbx)=\vep^{-n}\vp(\bbx/\vep),$
$\mathbbm{1}^{\vep}_{E^{\rM}}(\bbx)=\mathbbm{1}_{E^{\rM}}*\Phi_{\vep}(\bbx)$
and
$$N^K_{\rM}(t)=\sum_{\bbm\in\BZ^n}\mathbbm{1}^{\vep}_{E^{\rM}}\biggl(\f{\bbm-\bap}{t}\biggr)\,\,
\text{ for $\vep=\f{1}{tK}$, }$$ where $K\ge 1$ is a constant
independent of $t$ to be determined later.

One can see that the number of the integral lattice points inside an
ellipsoid is the same as that of the dilated lattice points (related
with the ellipsoid) inside a sphere. That is to say, since $\fL$ is
the lattice spanned by $\{\rT(e_1),\cdots,\rT(e_n)\}$ for the
standard basis $\{e_1,\cdots,e_n\}$ of $\BZ^n$, one has that
$$N^K_{\rM}(t)=\sum_{\bbm\in\BZ^n}\mathbbm{1}^{\vep}_{E^{\rM}}\biggl(\f{\bbm-\bap}{t}\biggr)
=\sum_{\bbm'\in\fL}\mathbbm{1}^{\vep}_{B}\biggl(\f{\bbm'-\bap'}{t}\biggr)\,\,
\text{ for $\vep=\f{1}{tK}$, }$$ where $K\ge 1$ is a constant
independent of $t$ to be determined later and $\bap'=\rT(\bap)$ is
the image of $\bap$ under the linear transformation $\rT$ on $\BR^n$
given by $\sqrt\rM$.

We recall the Poisson summation formula
\begin{equation}\sum_{\bbm'\in\fL}g(\bbm')=\f{1}{|\Pi_{\fL}|}\sum_{\bbm^*\in\fL^*}\widehat g(\bbm^*)
\,\,\text{ for $g\in\cS(\BR^n)$,}
\end{equation}
where $|\Pi_{\fL}|$ is the volume of the fundamental domain
$\Pi_{\fL}$ of the lattice $\fL$. Here we observe that
$\bbm^*=\rT^*(\bbm)$ is the image of $\bbm$ under the linear
transformation $\rT^*$, and so the dual lattice $\fL^*$ is spanned
by $\{\rT^*(e_1),\cdots,\rT^*(e_n)\}$. Since $\rT^*=\rT^{-1}$, it
easily follows from this formula (2.2) that
\begin{equation}\begin{split}N^K_{\rM}(t)&=|E^{\rM}|\,t^n+\f{t^n}{\det\rT}\sum_{\bbm^*\in\fL^*\s\{\bold 0\}}e^{-2\pi
i\bap'\cdot\bbm^*}\,\widehat{\mathbbm{1}_{B}}(t\bbm^*)\,\widehat\vp\biggl(\f{\bbm^*}{K}\biggr)\\
&=|E^{\rM}|\,t^n+\f{t^n}{\det\rT}\sum_{\bbm\in\BZ^n\s\{\bold
0\}}e^{-2\pi
i\bap\cdot\bbm}\,\widehat{\mathbbm{1}_{B}}(t\,\rT^*(\bbm))\,\widehat\vp\biggl(\f{\rT^*(\bbm)}{K}\biggr).
\end{split}\end{equation}
The Fourier coefficients of $\mathbbm{1}_{B}$ can be obtained via
the Bessel function as follows;
$$\widehat{\mathbbm{1}_{B}}(\bxi)=|\bxi|^{-\f{n}{2}}J_{\f{n}{2}}(2\pi|\bxi|)=\sum_{\el=0}^N
P_{\el}\,|\bxi|^{-\f{n+1}{2}-\el}\cos(2\pi|\bxi|+\phi_{\el})+\cO(|\bxi|^{-\f{n+1}{2}-N-1}),$$
\begin{equation}P_0=\f{1}{\pi}\,\,\,\text{ and
}\,\,\,\phi_0=-\f{n+1}{4}\,\pi,\end{equation} where $P_{\el}$ and
$\phi_{\el}, \el=1,2,\cdots,$ are real coefficients and phases,
respectively. The above asymptotics are valid for all $N\in\BN$. We
define the regularization of the error function by
$$F^K_{\rM}(t)=\f{N^K_{\rM}(t)-|E^{\rM}|\,t^n}{t^{\f{n-1}{2}}}.$$
By (2.3) and (2.4), we express $F^K_{\rM}(t)$ as the sum
$F^K_{\rM}(t)=\sum_{\el=0}^N\underline F_{\rM}^{K,\el}(t)+\underline
R_{\rM}^{K,N}(t),$ where
$$\underline F_{\rM}^{K,\el}(t)=\f{P_{\el}t^{-\el}}{\det\rT}\sum_{\bbm\in\BZ^n\s\{\bold
0\}}\f{\cos(2\pi
t|\rT^*(\bbm)|+\phi_{\el})}{|\rT^*(\bbm)|^{\f{n+1}{2}+\el}}\,e^{-2\pi
i\bap\cdot\bbm}\,\widehat\vp\biggl(\f{\rT^*(\bbm)}{K}\biggr)\,\,\text{
for $\el\ge 0$}.$$ With the condition $N>(n-1)/2$, the function
$\underline R_{\rM}^{K,N}(t)$ is continuous in $t>0$ and it
satisfies the bound
\begin{equation}|\underline R_{\rM}^{K,N}(t)|\ll t^{-N}\,\,\text{ for $N>\f{n-1}{2}$,
}\end{equation} uniformly in the parameter $K$. For convenience, we
truncate the sums $\underline F_{\rM}^{K,\el}$ by splitting it into
a sum over $|\bbm|_{\widehat\rM}\le K^{1+\zt/2}$ and a sum over
$|\bbm|_{\widehat\rM}>K^{1+\zt/2}$ with any $\zt>0$. Since
$|\widehat\vp(\bxi)|\ll(1+|\bxi|)^{-H}$ for an arbitrary $H>0$, the
sum over $|\bbm|_{\widehat\rM}>K^{1+\zt/2}$ is bounded by
\begin{equation}K^H\sum_{|\bbm|_{\widehat\rM}>K^{1+\zt/2}}|\bbm|_{\widehat\rM}^{-\f{n+1}{2}-\el-H}\ll
K^{-\f{\zt H}{2}+n}.\end{equation} For $\el\ge 0$, we set
\begin{equation}F_{\rM}^{K,\el}(t)=\f{P_{\el}t^{-\el}}{\det\rT}\sum_{\substack{\bbm\in\BZ^n\s\{\bold
0\}\\|\bbm|_{\widehat\rM}\le K^{1+\zt/2}}}\f{\cos(2\pi
t|\rT^*(\bbm)|+\phi_{\el})}{|\rT^*(\bbm)|^{\f{n+1}{2}+\el}}\,e^{-2\pi
i\bap\cdot\bbm}\,\widehat\vp\biggl(\f{\rT^*(\bbm)}{K}\biggr)\end{equation}
and we denote by $R_{\rM}^{K,N}(t)$ the remaining part of the sum.
Thus by (2.5) and (2.6) one sees that
\begin{equation}F^K_{\rM}(t)=\sum_{\el=0}^N F_{\rM}^{K,\el}(t)+R_{\rM}^{K,N}(t),
\end{equation} where
\begin{equation}|R_{\rM}^{K,N}(t)|\ll t^{-N}+K^{-H}\,\,\text{ for $t>0$ and
$N>\f{n-1}{2}$,}
\end{equation} with an arbitrary $H>0$. Since the Fourier transform of a radial function is radial,
it follows from (1.3), Lemma 2.1 and (2.7) that
\begin{equation}\begin{split}F^{K,\el}_{\rM}(t)&=\f{P_{\el}(\det\rM)^{\f{n-1+2\el}{4}}}{t^{\el}}\\
&\qquad\times\sum_{p=1}^{K^{2+\zt}}\,\f{\cos\biggl(\displaystyle\f{2\pi
t\sqrt
p}{(\det\rM)^{\f{1}{2}}}+\phi_{\el}\biggr)}{p^{\f{n+1+2\el}{4}}}\,\widehat
\vp\biggl(\f{\sqrt
p}{K(\det\rM)^{\f{1}{2}}}\biggr)\,\fr[\widehat\rM,\bap](p).
\end{split}\end{equation}

\section{Mean square value of exponential sums, extending work of
Marklof ( Proof of Theorem 1.3 ) }

In this section, we extend the work of Marklof \cite{M3} on the
Euclidean balls to some ellipsoids and prove Theorem 1.3 by applying
it.

We recall Shale-Weil representation and Theta sums. Most of
statements in this section  are elliptic adaptations of the results
taken from \cite{LV}, \cite{M1} and \cite{M2}. But we are dealing
with a more general situation than that in \cite{M1} and \cite{M2},
and thus we provide the proofs which are proper to reach our
results, even though the statements are similar to those in
\cite{M1} and \cite{M2}.

\subsection{Schr\"odinger representation and Shale-Weil representation}

Let $\om$ be the standard symplectic form on $\BR^{2n}$, i.e.
$$ \om\biggl(\begin{pmatrix} \bbx \\ \bby \end{pmatrix},\begin{pmatrix} \bbx' \\ \bby' \end{pmatrix}\biggr)
= \bbx \cdot \bby' - \bby \cdot \bbx'.$$ The Heisenberg group
$\BH(\BR^n)$ is defined to be the group $\BR^{2n} \ltimes \BR$ with
multiplication law:
$$ (\bxi, t)(\bxi', t') = (\bxi + \bxi' , t + t' + \om(\bxi, \bxi')/2).$$
We can express an element of the Heisenberg group as
$$ \biggl(\begin{pmatrix} \bbx \\ \bby \end{pmatrix}, t\biggr) = \biggl(\begin{pmatrix} \bbx \\
 \bold 0 \end{pmatrix}, 0\biggr) \biggl(\begin{pmatrix} \bold 0 \\ \bby \end{pmatrix}, 0\biggr)
 \biggl(\begin{pmatrix}\bold 0
 \\ \bold 0 \end{pmatrix}, t- \frac{\bbx\cdot\bby}{2}\biggr). $$
The Schr\"odinger representation of $\BH(\BR^n)$ on $f \in
L^2(\BR^n)$ is given by
\begin{equation}\rW\biggl( \begin{pmatrix} \bbx \\
\bby\end{pmatrix}, t\biggr) f (\mathbf w ) = \fe
\biggl(-t+\frac{\bbx \cdot \bby}{2}\biggr) \fe (-\bbx \cdot \mathbf
w) f(\mathbf w - \bby).\end{equation}

Let ${\rSp}_n(\BR)$ be the symplectic group of $2n \times 2n$
matrices which preserves $\om$, i.e.
$${\rSp}_n(\BR) = \left\{ \begin{pmatrix} \rA & \rB \\ \rC & \rD \end{pmatrix}\in{\rm GL}_{2n}(\BR)
:  \rC^t \rA = \rA^t \rC, \rA^t\rD-\rC^t \rB = {\rm I}_n, \rD^t \rB
= \rB^t \rD\right\}.$$ For any element $\rS\in\rSp_n(\BR)$, one can
define a new representation $\rW_\rS$ of $\BH(\BR^n)$ by
$$ \rW_\rS(\bxi, t ) = \rW(\rS \bxi, t).$$
By Stone-von Neumann theorem, all such representations are unitarily
equivalent, i.e. there exists a unitary operator $\rR(\rS)$ such
that
$$ \rR(\rS) \rW(\bxi ,t) \rR(\rS)^{-1} = \rW(\rS\bxi ,t).$$
Then $\rS \mapsto \rR(\rS)$ is determined up to a unitary phase
factor, thus defines a projective representation, called the
projective {\it Shale-Weil representation} of the symplectic group,
i.e. $ \rR(\rS\rS') = c(\rS,\rS') \rR(\rS)\rR(\rS'),$ with cocycle
$c(\rS,\rS') \in \mathbb C$ with $|c(\rS,\rS') | = 1$.

Now consider a family $$\cM= \left\{\rS=\begin{pmatrix} \rA
 & \rB \\ \rC & \rD \end{pmatrix} \in\rSp_n(\BR)  :  \rA,\rB,\rC,\rD\in\fS_n(\BR) \right\}$$
where $\fS_n(\BR)$ is the family of all symmetric $n\times n$
matrices with real components. Then any element in $\cM$ satisfies
\begin{equation}\rC\rA=\rA\rC,\,\,\rA\rD-\rC\rB=\rI_n,\,\,\rD\rB=\rB\rD.\end{equation}
For $\rS\in\cM$, we shall give an explicit description of cocycles
as well as that of $\rR(\rS)$ for these elements. If $\rS_i =
\begin{pmatrix} \rA_i & \rB_i
\\ \rC_i & \rD_i \end{pmatrix}$ for $i=1,2,3$ with $\rS_1 \rS_2 =\rS_3$, then
\begin{equation}c(\rS_1, \rS_2 ) =\mathrm{exp}\biggl(\frac{-\pi i\,\tau(\BR_{\bbx}^n,
\rS_1 \BR_{\bbx}^n, \rS_3 \BR_{\bbx}^n)}{4}\biggr),\end{equation}
where $\BR_{\bbx}^n=\ds\biggl\{\begin{pmatrix} \bbx \\ \bold 0
\end{pmatrix}\in\BR^{2n}:\bbx\in\BR^n\biggr\}$ and $\tau(\BR_{\bbx}^n, \rS_1 \BR_{\bbx}^n, \rS_3
\BR_{\bbx}^n )$ is the Maslov index (see \cite{LV}). In particular,
if $\rC_1,\rC_2$ are invertible, the corresponding cocycle is
\begin{equation}c(\rS_1,\rS_2) =
\mathrm{exp} \left( \frac{-\pi i \,\,\sgn(\rC_1^{-1} \rC_3
\rC_2^{-1})}{4} \right).\end{equation} We also have the explicit
representation
$$ \rR(\rS) f (\bby) = \int_{\BR^n/\ker(\rC)}e^{\pi i\, \om(\rA\rB\bby,
\bby)} e^{-2 \pi i\,\om(\rB\bby,\rC\bbx)} e^{-\pi i
\,\om(\rD\rC\bbx, \bbx)}f(\rA\bby-\rC\bbx)\,d\bbx.$$ In particular,
it is represented as \begin{equation}\begin{split}&\rR(\rS) f
(\bbw)\\&=\left\{\begin{array}{cl}
     |\det \rA|^{\f{1}{2}}\,e^{\pi i\,\om(\rA\rB\bbw,
\bbw)} f(\rA\bbw), & \rC=\rD_{\bold 0}, \\
     \displaystyle\int_{\BR^n}
     e^{\pi i(\om(\rC^{-1}\rD\bbw',\bbw')+\om(\rA\rC^{-1}\bbw,\bbw)-2\om(\rC^{-1}\bbw,\bbw'))}
     \f{f(\bbw')}{|\det\rC|^{\f{1}{2}}}\,d\bbw' , & \det\rC\neq 0.
       \end{array}\right.\end{split}\end{equation}

\subsection{Cocycles and Iwasawa decomposition of $\rSL_2(\BR)^n$}

Let $\bG^n$ be a semi-direct product of $\rSL_2(\BR)^n$ and
$\BR^{2n}$ induced from an action of $\rSL_2(\BR)^n$ on $\BR^{2n}$
as follows; for $\rM=\biggl(\begin{pmatrix} a_1 & b_1\\ c_1 &
d_1\end{pmatrix},\cdots,\begin{pmatrix} a_n & b_n\\ c_n &
d_n\end{pmatrix}\biggr)\in\rSL_2(\BR)^n$ and $\bxi=\begin{pmatrix}
\bbx\\ \bby\end{pmatrix}\in\BR^{2n}$, we define
$$\rM\,\bxi=\begin{pmatrix} \rD_{\bba} & \rD_{\bbb}\\ \rD_{\bbc} &
\rD_{\bbd}\end{pmatrix}\begin{pmatrix} \bbx \\ \bby \end{pmatrix}
$$ where
$\bba=(a_1,\cdots,a_n),\bbb=(b_1,\cdots,b_n),\bbc=(c_1,\cdots,c_n),
\bbd=(d_1,\cdots,d_n)\in\BR^n$. We note that
$\bG^n=\rSL_2(\BR)^n\ltimes\BR^{2n}\cong[\rSL_2(\BR)\ltimes\BR^2]^n$
has the following multiplication law
$(\rM;\bxi)(\rM';\bxi')=(\rM\rM';\bxi+\rM\bxi').$

\begin{prop} If $\rM=\begin{pmatrix} \rD_{\bba} & \rD_{\bbb} \\
\rD_{\bbc} & \rD_{\bbd} \end{pmatrix}$, $\rM'=\begin{pmatrix}
\rD_{\bba'} & \rD_{\bbb'} \\ \rD_{\bbc'}
& \rD_{\bbd'} \end{pmatrix}$, $\rM''=\begin{pmatrix} \rD_{\bba''} & \rD_{\bbb''} \\
\rD_{\bbc''} & \rD_{\bbd''} \end{pmatrix}\in\rSL_2(\BR)^n$ with
$\rM\rM'=\rM''$, then the corresponding cocycle is
$$c(\rM,\rM')=\exp\biggl(-\f{i\pi}{4}\sum_{k=1}^n\sgn(c_k c'_k
c''_k)\biggr)$$ where $\sgn(x)=x/|x|$ for $x\in\BR\s\{0\}$ and
$\sgn(0)=0$.
\end{prop}

\pf We note that
$c(\rM,\rM')=\exp(-\f{i\pi}{4}\,\tau(\BR^n_{\bbx},\rM\BR^n_{\bbx},\rM\rM'\BR^n_{\bbx}))$
and we recall that the Maslov index $\tau(V_1,V_2,V_3)$ is the
signature of the quadratic form $Q(\bbx_1+\bbx_2+\bbx_3)$ on
$V_1\oplus V_2\oplus V_3$ defined by
$$Q(\bbx_1+\bbx_2+\bbx_3)=\om(\bbx_1,\bbx_2)+\om(\bbx_2,\bbx_3)+\om(\bbx_3,\bbx_1)$$
where $V_1,V_2,V_3$ are three Lagrangian planes in a
$2n$-dimensional symplectic vector space $V$(see \cite{LV}). Since
$\bG^n$ is the direct product of $\bG^1=\rSL_2(\BR)\ltimes\BR^2$ and
$\rM=\begin{pmatrix} \rD_{\bba} & \rD_{\bbb} \\
\rD_{\bbc} & \rD_{\bbd} \end{pmatrix}\in\rSL_2(\BR)^n$ acts on
$\BR^n_{\bbx}\subset\BR^{2n}$ componentwise, we may decompose
$\BR^n_{\bbx}$ and $\BR^{2n}$ as follows;
$$\BR^n_{\bbx}=\BR_{x_1}\oplus\cdots\oplus\BR_{x_n}\subset\BR^{2n}=\BR^2\oplus\cdots\oplus\BR^2,$$
where $\BR_{x_k}=\ds\biggl\{\begin{pmatrix} x_k \\
0 \end{pmatrix}\in\BR^2: x_k\in\BR\biggr\}$ and
each $\begin{pmatrix} a_k & b_k \\
c_k & d_k \end{pmatrix}\in\rSL_2(\BR)$ acts on
$\BR_{x_k}\subset\BR^2$ for $k=1,\cdots,d$. Then the Maslov index
$\tau(\BR^n_{\bbx},\rM\BR^n_{\bbx},\rM\rM'\BR^n_{\bbx})$ is of the
form
$$\tau(\BR^n_{\bbx},\rM\BR^n_{\bbx},\rM\rM'\BR^n_{\bbx})=\sum_{k=1}^n
\tau\biggl(\BR_{x_k},\begin{pmatrix} a_k & b_k\\ c_k &
d_k\end{pmatrix}\BR_{x_k},\begin{pmatrix} a_k & b_k\\ c_k &
d_k\end{pmatrix}\begin{pmatrix} a'_k & b'_k\\ c'_k &
d'_k\end{pmatrix}\BR_{x_k}\biggr).$$ Thus it suffices to prove when
$\bG^n=\bG^1$. In case that both $c$ and $c'$ are not zero, it
follows from (3.4). If $c=0$, that is, $\rM=\begin{pmatrix} a & b\\
0 & d
\end{pmatrix}$, $\rM'=\begin{pmatrix} a' & b'\\ c' & d'
\end{pmatrix}$ $\rM=\begin{pmatrix} a'' & b''\\ c'' & d''
\end{pmatrix}$, then it is easy to check that $\rM\BR_x=\BR_x$, and
so $\om(\bbx_1,\bbx_2)\equiv 0$ and
$\sgn[\om(\bbx_2,\bbx_3)]=-\sgn[\om(\bbx_3,\bbx_1)]$, where
$\bbx_1\in\BR_x$, $\bbx_2\in\rM\BR_x$ and $\bbx_3\in\rM\rM'\BR_x$.
Similarly, if $c'=0$, then we easily get that
$\rM\BR_x=\rM\rM'\BR_x$, $\om(\bbx_2,\bbx_3)\equiv 0$ and
$\sgn[\om(\bbx_1,\bbx_2)]=-\sgn[\om(\bbx_3,\bbx_1)]$. Therefore the
required equality holds. \qed

Any $2n\times 2n$ matrix $\rM=\begin{pmatrix} \rD_\bba & \rD_\bbb \\
\rD_\bbc & \rD_\bbd \end{pmatrix}\in\rSL_2(\BR)^n$ admits the unique
Iwasawa decomposition
\begin{equation}\rM=\begin{pmatrix} \rI_n & \rD_\bbu \\
\rD_{\bold 0} & \rI_n \end{pmatrix}\begin{pmatrix} \rD_{\bbv^{1/2}}  & \rD_{\bold 0} \\
\rD_{\bold 0} & \rD_{\bbv^{-1/2}} \end{pmatrix}\begin{pmatrix} \rD_{\cos\bphi} & -\rD_{\sin\bphi} \\
\rD_{\sin\bphi} & \rD_{\cos\bphi}
\end{pmatrix}:=(\bbz,\bphi),\end{equation} where $\bold
0=(0,\cdots,0)\in\BR^n$, $\bold 1=(1,\cdots,1)\in\BR^n$,
$\bbu=(u_1,\cdots,u_n)$, $\bbv^{1/2}=(v_1^{1/2},\cdots,v_n^{1/2})$,
$\bbv^{-1/2}=(v_1^{-1/2},\cdots,v_n^{-1/2})$,
$\cos\bphi=(\cos\phi_1,\cdots,\cos\phi_n)$,
$\sin\bphi=(\sin\phi_1,\cdots,\sin\phi_n)$, $\bbz=(z_1,\cdots,z_n)$,
$\bphi=(\phi_1,\cdots,\phi_n)$, $\phi_k\in[0,2\pi)$ for
$k=1,\cdots,n$ and $z_k=u_k+i v_k\in\BH^2:=\{z\in\BC:\Ima(z)>0\}$.

\subsection{Jacobi's theta sums}

The Jacobi group is defined as the semidirect product $\rSp_n(\BR)
\ltimes \BH(\BR^n),$ with multiplication law
$$ (\rS;\bxi,t)(\rS';\bxi',t') = (\rS\rS';\bxi+\rS\bxi',t+ t'+\frac{1}{2}\om(\bxi,\rS\bxi')).$$
Then $\rR(\rS;\bxi,t) := \rW(\bxi,t)\rR(\rS)$ defines a projective
representation of the Jacobi group, called {\it the
Schr\"odinger-Weil representation.}

For $f \in \cS(\BR^n)$, we note that
$$\lim_{\phi_1 \to 0^{\pm}}\cdots\lim_{\phi_n \to 0^{\pm}} \rR(i\bold
1,\bphi)f(\mathbf w) =e^{\pm\f{i n\pi}{4}} f(\mathbf w),$$ and thus
this projective representation is not continuous at
$\bphi=(\phi_1,\cdots,\phi_n)\in \BZ^n \pi$. This obstacle can be
overcome by putting $$\widetilde \rR (\bbz, \bphi)=e^{-\f{i
\pi}{4}(\sm_{\phi_1}+\cdots+\sm_{\phi_n})}\,\rR(\bbz, \bphi),$$
where $\bbz=(z_1,\cdots,z_n)\in[\BH^2]^n$, and $\sigma_{\phi_k} =
2\nu$ if $\phi_k=\nu\pi$ and $\sigma_{\phi_k}= 2\nu+1$ if $\nu\pi <
\phi_k < (\nu+1)\pi$ (see \cite{M1}). Write
$\widetilde\rR(\bbz,\bphi\,;\bxi,t)=\rW(\bxi,t)\widetilde
\rR(\bbz,\bphi),$ so that one has that
$$\widetilde\rR(\bbz,\bphi\,;\bxi,t)=e^{-\f{i\pi}{4}(\sm_{\phi_1}+\cdots+\sm_{\phi_n})}\,\rW(\bxi,
t)\,\rR(\bbz, \bphi).$$ Now we define the Jacobi's theta sum for
$f\in\cS (\BR^n)$ by
$$\Theta_f(\bbz, \bphi\,; \bxi, t) = \sum_{\bbm \in \BZ^n} [\,\widetilde\rR(\bbz, \bphi\,; \bxi, t) f\,](\bbm).$$
More explicitly, for $(\bbz,\bphi)\in\rSL_2(\BR)^n$ and $\bxi =
\begin{pmatrix} \bbx \\ \bby \end{pmatrix}\in\BR^{2n},$ by (3.1) we obtain that
\begin{equation}\begin{split}
&\Theta_f(\bbz,\bphi\,; \bxi, t)= (v_1\cdots v_n)^{\f{1}{4}}
\,\fe\bigl(-t+\frac{\bbx\cdot\bby}{2}\bigr)\\
&\times\sum_{ \bbm \in \BZ^n}
f_{\bphi}((m_1-y_1)v_1^{\f{1}{2}},\cdots,(m_n-y_n)v_n^{\f{1}{2}})
\,\fe\biggl(-\f{1}{2}\sum_{k=1}^n(m_k-y_k)^2
u_k-\bbm\cdot\bbx\biggr),
\end{split}\end{equation}
where $f_{\bphi}=\widetilde\rR(i\bold 1,\bphi)f$. It is easy to see
that if $f \in \cS(\BR^n)$, then $f_{\bphi} \in \cS(\BR^n)$ for
fixed $\bphi$, and thus $\widetilde{\rR}(\bbz, \bphi\,; \bxi ,t ) f
\in \cS(\BR^n)$ as well for fixed $(\bbz, \bphi\,; \bxi, t)$. Also
we note that $f_{\bold 0}=f$. This guarantees rapid convergence of
the above series.

The relevant discrete subgroup $\Gm^n$ of $\bG^n$ is defined by
$$\Gm^n=\biggl\{\biggl(\rM;\f{1}{2}\begin{pmatrix}\bba\bbb \\
\bbc\bbd\end{pmatrix}+\begin{pmatrix} \bbm \\
\bbh
\end{pmatrix}\biggr):\bba,\bbb,\bbc,\bbd,\bbm,\bbh\in\BZ^n,\rM=\begin{pmatrix} \rD_{\bba} & \rD_{\bbb} \\
\rD_{\bbc} & \rD_{\bbd} \end{pmatrix}\in\rSL_2(\BZ)^n\biggr\}
$$ where $\bba\bbb=(a_1 b_1,\cdots,a_n b_n)$ and $\bbc\bbd=(c_1
d_1,\cdots,c_n d_n)$. Then the left action of $\Gm^n$ on $\bG^n$ is
properly discontinuous. Moreover, a fundamental domain of $\Gm^n$ is
given by
$$\cF_{\Gm^n}=\prod_{k=1}^n\biggl(\cF_{\rSL_2(\BZ)}\times\{\phi_k\in[0,\pi)\}\times\biggl\{\begin{pmatrix} x_k \\
y_k\end{pmatrix}\in[-1/2,1/2)^2\biggr\}\biggr)$$ where
$\cF_{\rSL_2(\BZ)}=\{z=u+iv\in\BH^2:u\in[-\f{1}{2},\f{1}{2}),
|z|>1\}$ is the fundamental domain in $\BH^2$ of the modular group
$\rSL_2(\BZ)$.

As a matter of convenience, we write $\rM=\begin{pmatrix} \rD_{\bba} & \rD_{\bbb} \\
\rD_{\bbc} & \rD_{\bbd} \end{pmatrix}\in\rSL_2(\BR)^n$ as
$$\rM=\biggl(\begin{pmatrix} a_1 & b_1 \\
c_1 & d_1
\end{pmatrix},\cdots,\begin{pmatrix} a_n & b_n \\
c_n & d_n
\end{pmatrix}\biggr).$$ Then we observe that $\Gm^n$ is generated by
three types of matrices as follows;
$$\rM_{\bbh}=\biggl(\begin{pmatrix} 1 & 0 \\
0 & 1
\end{pmatrix},\begin{pmatrix} 1 & 0 \\
0 & 1
\end{pmatrix},\cdots,\begin{pmatrix} 1 & 0 \\
0 & 1
\end{pmatrix};\begin{pmatrix} \bbh_1 \\
\bbh_2
\end{pmatrix}\biggr)\,\,\text{ for $\bbh_1,\bbh_2\in\BZ^n$, }$$
\begin{equation*}\begin{split}\rF_1=&\biggl(\begin{pmatrix} 0 & -1 \\
1 & 0
\end{pmatrix},\begin{pmatrix} 1 & 0 \\
0 & 1
\end{pmatrix},\begin{pmatrix} 1 & 0 \\
0 & 1
\end{pmatrix},\cdots,\begin{pmatrix} 1 & 0 \\
0 & 1
\end{pmatrix};\begin{pmatrix} {\bold 0} \\
{\bold 0}
\end{pmatrix}\biggr),\\
\rF_2=&\biggl(\begin{pmatrix} 1 & 0 \\
0 & 1
\end{pmatrix},\begin{pmatrix} 0 & -1 \\
1 & 0
\end{pmatrix},\begin{pmatrix} 1 & 0 \\
0 & 1
\end{pmatrix},\cdots,\begin{pmatrix} 1 & 0 \\
0 & 1
\end{pmatrix};\begin{pmatrix} {\bold 0} \\
{\bold 0}
\end{pmatrix}\biggr),\\
&\cdots\\
\rF_n=&\biggl(\begin{pmatrix} 1 & 0 \\
0 & 1
\end{pmatrix},\begin{pmatrix} 1 & 0 \\
0 & 1
\end{pmatrix},\cdots,\begin{pmatrix} 1 & 0 \\
0 & 1
\end{pmatrix},\begin{pmatrix} 0 & -1 \\
1 & 0
\end{pmatrix};\begin{pmatrix} {\bold 0} \\
{\bold 0}
\end{pmatrix}\biggr),
\end{split}\end{equation*}
\begin{equation*}\begin{split}\rU_1=&\biggl(\begin{pmatrix} 1 & 1 \\
0 & 1
\end{pmatrix},\begin{pmatrix} 1 & 0 \\
0 & 1
\end{pmatrix},\begin{pmatrix} 1 & 0 \\
0 & 1
\end{pmatrix},\cdots,\begin{pmatrix} 1 & 0 \\
0 & 1
\end{pmatrix};\begin{pmatrix} \bbs_1 \\
{\bold 0}
\end{pmatrix}\biggr),\\
\rU_2=&\biggl(\begin{pmatrix} 1 & 0 \\
0 & 1
\end{pmatrix},\begin{pmatrix} 1 & 1 \\
0 & 1
\end{pmatrix},\begin{pmatrix} 1 & 0 \\
0 & 1
\end{pmatrix},\cdots,\begin{pmatrix} 1 & 0 \\
0 & 1
\end{pmatrix};\begin{pmatrix} \bbs_2 \\
{\bold 0}
\end{pmatrix}\biggr),\\
&\cdots\\
\rU_n=&\biggl(\begin{pmatrix} 1 & 0 \\
0 & 1
\end{pmatrix},\begin{pmatrix} 1 & 0 \\
0 & 1
\end{pmatrix},\cdots,\begin{pmatrix} 1 & 0 \\
0 & 1
\end{pmatrix},\begin{pmatrix} 1 & 1 \\
0 & 1
\end{pmatrix};\begin{pmatrix} \bbs_n \\
{\bold 0}
\end{pmatrix}\biggr),
\end{split}\end{equation*} where $\bbs_k$ is an element of $\BR^n$
whose entries are all $0$ except that the $k$-th entry is $1/2$.

\begin{lemma} For $f\in\cS(\BR^n)$, the Jacobi's theta sum $\Theta_f (\bbz, \bphi\,; \bxi,
t)$ satisfies the following properties.

$(a)$ $\Theta_f(\rF_k(\bbz, \bphi\,; \bxi,
t))=e^{-\f{i\pi}{4}}\Theta_f(\bbz, \bphi\,; \bxi, t)$ for
$k=1,\cdots,n,$

$(b)$ $\Theta_f(\rU_k(\bbz, \bphi\,; \bxi,
t))=e^{-\f{i\pi}{2}y_k}\Theta_f (\bbz, \bphi\,; \bxi, t)$ for
$k=1,\cdots,n$,

$(c)$ $\Theta_f(\rM_{\bbh}(\bbz, \bphi\,; \bxi,
t))=e^{-i\pi(\bbh_1\cdot\bbh_2)}\Theta_f (\bbz, \bphi\,; \bxi, t)$
for $\bbh_1,\bbh_2\in\BZ^n$ and $\bxi =
\begin{pmatrix} \bbx \\ \bby \end{pmatrix}.$
\end{lemma}

\pf Take any $f\in\cS(\BR^n)$ and $k=1,\cdots,n$.

(a) By (3.5), we easily see that
\begin{equation*}[\rR(\rF_k)f](x_1,\cdots,x_{k-1},m_k,x_{k+1},\cdots,x_n)
=\int_{\BR}e^{-2\pi i m_k x_k}f(x_1,\cdots,x_n)\,dx_k,
\end{equation*} which is the partial Fourier transformation of $f$
in terms of $x_k$-variable. Since $\widetilde\rR(\bbz, \bphi\,;
\bxi, t)f\in\cS(\BR^n)$ for fixed $(\bbz, \bphi\,; \bxi, t)$, it
follows from the Poisson summation formula that
\begin{equation}\begin{split}
\sum_{\bbm\in\BZ^n}[\rR(\rF_k)\widetilde\rR(\bbz, \bphi\,; \bxi,
t)f](\bbm)&=\sum_{\bbm'\in\BZ^{n-1}}\sum_{m_k\in\BZ}[\rR(\rF_k)\widetilde\rR(\bbz,
\bphi\,; \bxi, t)f](\bbm)\\
&=\sum_{\bbm'\in\BZ^{n-1}}\sum_{m_k\in\BZ}[\widetilde\rR(\bbz,
\bphi\,; \bxi, t)f](\bbm)\\
&=\sum_{\bbm\in\BZ^n}[\widetilde\rR(\bbz, \bphi\,; \bxi, t)f](\bbm),
\end{split}\end{equation}
where $\bbm'=(m_1,\cdots,m_{k-1},m_{k+1},\cdots,m_n)\in\BZ^{n-1}$.

On the other hand, using Proposition 3.3 and the fact that
$0<\theta_k=\arg(z_k)<\pi$ for $z_k\in\BH^2$, one obtains that
\begin{equation*}\begin{split}
&\rR(\rF_k)\widetilde\rR(\bbz, \bphi\,; \bxi,
t)\\&=\rR(\rF_k)\rW(\bxi,t)\rR(\bbz,\bold 0)\widetilde\rR(i\bold 1,
\bphi)=\rW(\rF_k\bxi,t)\rR(\rF_k)\rR(\bbz,\bold
0)\widetilde\rR(i\bold 1,
\bphi)\\
&=\rW(\rF_k\bxi,t)\rR((z_1,0),\cdots,(z_{k-1},0),(\f{-1}{z_k},\theta_k),(z_{k+1},0),\cdots,(z_n,0))\widetilde\rR(i\bold
1, \bphi)\\
&=\rW(\rF_k\bxi,t)\rR((z_1,0),\cdots,(z_{k-1},0),(\f{-1}{z_k},0),(z_{k+1},0),\cdots,(z_n,0))\\
&\qquad\qquad\quad\times\rR((i,0),\cdots,(i,0),(i,\theta_k),(i,0),\cdots,(i,0))\,\,\widetilde\rR(i\bold
1, \bphi)\\
&=\rW(\rF_k\bxi,t)\rR((z_1,0),\cdots,(z_{k-1},0),(\f{-1}{z_k},0),(z_{k+1},0),\cdots,(z_n,0))\\
&\qquad\qquad\quad\times
e^{i\f{\pi}{4}}\,\widetilde\rR((i,\phi_1),\cdots,(i,\phi_{k-1}),(i,\phi_k+\theta_k),(i,\phi_{k+1}),\cdots,(i,\phi_n))\\
&=e^{i\f{\pi}{4}}\,\widetilde\rR((z_1,\phi_1),\cdots,(z_{k-1},\phi_{k-1}),(\f{-1}{z_k},\phi_k+\theta_k),(z_{k+1},\phi_{k+1}),\cdots,(z_n,\phi_n);
\rF_k\bxi,t)\\
&=e^{i\f{\pi}{4}}\,\widetilde\rR(\rF_k(\bbz, \bphi\,; \bxi, t)).
\end{split}\end{equation*}
Hence the required equality (a) follows from (3.7).

(b) A simple computation for $\bbm\in\BZ^n$ yields that
\begin{equation*}\begin{split}
[\rR(\rU_k)f](\bbm)&=\biggl[\rW\biggl(\begin{pmatrix} \bbs_k \\
\bold 0
\end{pmatrix},0\biggr)\widetilde\rR((i,0),\cdots,(i,0),(i+1,0),(i,0),\cdots,(i,0))f\biggr](\bbm)\\
&=\fe(-\bbs_k\cdot\bbm)\fe\biggl(-\f{1}{2}m_k^2\biggr)f(\bbm)=\fe\biggl(-\f{m_k(m_k+1)}{2}\biggr)f(\bbm)=f(\bbm),
\end{split}\end{equation*}
since $m_k(m_k+1)\in 2\BZ$. Thus by replacing $f$ by
$\widetilde\rR(\bbz, \bphi\,; \bxi, t)f$, we get that
\begin{equation}\sum_{\bbm\in\BZ^n}[\widetilde\rR(\rU_k)\widetilde\rR(\bbz,
\bphi\,; \bxi, t)f](\bbm) =\sum_{\bbm\in\BZ^n}[\widetilde\rR(\bbz,
\bphi\,; \bxi, t)f](\bbm).\end{equation} Since $[\widetilde\rR(\bbz,
\bphi\,; \bxi, t+s)f](\bbm)=e^{2i\pi s}[\widetilde\rR(\bbz, \bphi\,;
\bxi, t)f](\bbm)$ and
\begin{equation*}\begin{split}
&\widetilde\rR(\rU_k)\widetilde\rR(\bbz, \bphi\,; \bxi,
t)\\&=\widetilde\rR((z_1,\phi_1),\cdots,(z_{k-1},\phi_{k-1}),(z_k+1,\phi_k),(z_{k+1},\phi_{k+1}),\cdots,(z_n,\phi_n);\\
&\biggl(\begin{pmatrix} x_1 \\
y_1
\end{pmatrix},\cdots,\begin{pmatrix} x_{k-1} \\
y_{k-1}\end{pmatrix},\begin{pmatrix} \f{1}{2} \\
0
\end{pmatrix}+\begin{pmatrix} 1 & 1 \\
0 & 1
\end{pmatrix}\begin{pmatrix} x_k \\
y_k
\end{pmatrix},\begin{pmatrix} x_{k+1} \\
y_{k+1}\end{pmatrix},\cdots,\begin{pmatrix} x_n \\
y_n\end{pmatrix}\biggr),t+\f{1}{4}y_k),
\end{split}\end{equation*} we obtain that
\begin{equation}[\widetilde\rR(\rU_k)\widetilde\rR(\bbz, \bphi\,; \bxi,
t)f](\bbm)=e^{\f{i\pi}{2}y_k}\widetilde\rR(\rU_k(\bbz, \bphi\,;
\bxi, t))f](\bbm).\end{equation} Thus we can complete the part (b)
by applying (3.8) and (3.9).

(c) If $\bbm\in\BZ^n$, then it follows from simple computation that
\begin{equation}\begin{split}
&[\widetilde\rR(\rM_{\bbn}(\bbz, \bphi\,; \bxi, t))f](\bbm)\\
&=[\widetilde\rR\biggl(\bbz,\bphi;\begin{pmatrix} \bbh_1 \\
 \bbh_2
\end{pmatrix}+\begin{pmatrix} \bbx \\
\bby
\end{pmatrix},t+\f{1}{2}(\bbn_1\cdot\bby-\bbn_2\cdot\bbx)\biggr)f](\bbm)\\
&=[\rW\biggl(\begin{pmatrix} \bbh_1 \\
\bbh_2
\end{pmatrix},0\biggr)\rW\biggl(\begin{pmatrix} \bbx \\
 \bby
\end{pmatrix},t\biggr)\widetilde\rR(\bbz, \bphi)f](\bbm)\\
&=[\rW\biggl(\begin{pmatrix} \bbh_1 \\
\bbh_2
\end{pmatrix},0\biggr)\widetilde\rR(\bbz, \bphi\,; \bxi,
t)f](\bbm)\\
&=e^{-i\pi(\bbh_1\cdot\bbh_2)}e^{2i\pi(\bbh_1\cdot\bbm)}[\widetilde\rR(\bbz,
\bphi\,; \bxi, t)f](\bbm-\bbh_2)\\
&=e^{-i\pi(\bbh_1\cdot\bbh_2)}[\widetilde\rR(\bbz, \bphi\,; \bxi,
t)f](\bbm-\bbh_2).
\end{split}\end{equation}
Hence the required equality (c) can be achieved by adding up (3.10)
on $\bbm\in\BZ^n$. \qed

\begin{prop} For $f,g\in\cS(\BR^n)$, $\Theta_f (\bbz, \bphi\,; \bxi, t) \overline{\Theta_g (\bbz,
\bphi\,;\bxi,t)}$ is invariant under the left action of $\Gm^n$.
\end{prop}
\pf It suffices to show that $\Theta_f (\bbz, \bphi\,; \bxi, t)
\overline{\Theta_g (\bbz, \bphi\,;\bxi,t)}$ is invariant under the
left action of generators of $\Gm^n$. Hence this immediately follows
from Lemma 3.5. \qed

We are now interested in the asymptotic limit of the mean square
value of the exponential sum $\fr[\rD_{\bba},\bap](p)$ given by
\begin{equation*}\fr[\rD_{\bba},\bap](p)=\sum_{\substack{\bbm\in\BZ^n \\
|\bbm|^2_{\bba}=p}} \fe(\bbm\cdot\bap)
\end{equation*} where $|\bbm|^2_{\bba}=Q_{\rD_{\bba}}(\bbm)$ for $\bba\in\BN^n$, $p\in\BN$,
$\bap\in\BR^n$ and $\fe(t)=e^{-2\pi i t}$ for $t>0$.

In order to obtain the following asymptotic limit of the mean square
value of the exponential sums
\begin{equation*}\lim_{N\to\iy}\f{1}{N^{n/2}}\rR[\rD_{\bba},\bap](N):=\lim_{N\to\iy}\f{1}{N^{n/2}}
\sum_{p=1}^N|\fr[\rD_{\bba},\bap](p)|^2,
\end{equation*}
we first consider the explicit form of $\Th_f\overline\Th_g$ on
$\rSL_2(\BR)^n\ltimes\BR^{2n}$ (when $f,g\in\cS(\BR^n)$) as follows;
\begin{equation*}\begin{split}
&\Theta_f(\bbz,\bphi\,; \bxi, t)\overline{\Theta_g(\bbz,\bphi\,;
\bxi, t)}=(v_1\cdots v_n)^{\f{1}{2}}\\
&\times\sum_{\bbm,\bbh \in \BZ^n}
f_{\bphi}((m_1-y_1)v_1^{\f{1}{2}},\cdots,(m_n-y_n)v_n^{\f{1}{2}})\overline{g_{\bphi}((h_1-y_1)v_1^{\f{1}{2}},
\cdots,(h_n-y_n)v_n^{\f{1}{2}})}\\
&\times\fe\biggl(-\f{1}{2}\sum_{k=1}^n(m_k-y_k)^2
u_k-\bbm\cdot\bbx\biggr)\fe\biggl(\f{1}{2}\sum_{k=1}^n(h_k-y_k)^2
u_k+\bbh\cdot\bbx\biggr),
\end{split}\end{equation*}
where $f_{\bphi}=\widetilde\rR(i\bold 1,\bphi)f$ and
$g_{\bphi}=\widetilde\rR(i\bold 1,\bphi)g$.

We set $\bbz=\bba z$ and $\bphi=\bold 0$ for $\bba\in\BN^n$ and
$z=u+i v\in\BH^2$, and $\bxi=\begin{pmatrix} \bbx \\  \bold 0
\end{pmatrix}.$ Choose some $f,g\in\cS(\BR^n)$ to be a
rapidly decreasing function $f(|\bbw|^2)=g(|\bbw|^2)=\psi(|\bbw|^2)$
which approximate the characteristic function
$\mathbbm{1}_{[0,1]}(|\bbw|^2)$ for $\bbw\in\BR^n$. Since $f_{\bold
0}=f$ and $g_{\bold 0}=g$, we have that
\begin{equation*}\begin{split}
&\bigl|\Theta_f(\bba z,\bold 0\,; \bxi,
t)\bigr|^2=(a_1\cdots a_n)^{\f{1}{2}}v^{\f{n}{2}}\\
&\qquad\qquad\times\sum_{\bbm,\bbh \in
\BZ^n}\mathbbm{1}_{[0,1]}(|\bbm|_{\bba}^2\,
v)\,\mathbbm{1}_{[0,1]}(|\bbh|_{\bba}^2\,v)\,\fe\biggl(\f{1}{2}\,u\bigl(|\bbh|^2_{\bba}-|\bbm|^2_{\bba}\bigr)
+\bigl(\bbh-\bbm\bigr)\cdot\bbx\biggr).
\end{split}\end{equation*}
Integrating $\bigl|\Theta_f(\bba z,\bold 0\,; \bxi, t)\bigr|^2$ on
$[0,2]$ in terms of $u$, one has that
\begin{equation}\begin{split}\f{1}{2}\int_0^2\bigl|\Theta_f(\bba z,\bold 0\,; \bxi,
t)\bigr|^2\,du=(a_1\cdots
a_n)^{\f{1}{2}}v^{\f{n}{2}}\sum_{\substack{\bbm,\bbh\in\BZ^n\\|\bbm|_{\bba}^2=|\bbh|_{\bba}^2\le
1/v}}\fe\bigl((\bbh-\bbm)\cdot\bbx\bigr). \end{split}\end{equation}
On the other hand, we have that
\begin{equation}|\fr[\rD_{\bba},\bap](p)|^2=\sum_{\substack{\bbm,\bbh\in\BZ^n\\|\bbm|_{\bba}^2=|\bbh|_{\bba}^2=p}}
\fe\bigl((\bbh-\bbm)\cdot\bap\bigr).
\end{equation}
If we set $N=[1/v]$ and $\bbx=\bap$ in (3.12) and (3.13), then we
obtain that
\begin{equation}\lim_{N\to\iy}\f{1}{N^{n/2}}\sum_{p=1}^N|\fr[\rD_{\bba},\bap](p)|^2=(a_1\cdots
a_n)^{-\f{1}{2}}\lim_{v\to 0}\f{1}{2}\int_0^2\bigl|\Theta_f(\bba
z,\bold 0\,; \bxi, t)\bigr|^2\,du.
\end{equation}

\subsection{Equidistribution of closed orbits}

Consider a subgroup
$$\biggl\{\rM\in\rSL_2(\BR)^n:\rM:=\begin{pmatrix} \rD_{\bba} & \rD_{\bbb} \\
\rD_{\bbc} & \rD_{\bbd}
\end{pmatrix}=\begin{pmatrix} a\,\rI_n & b\,\rI_n \\
c\,\rI_n & d\,\rI_n
\end{pmatrix}, ad-bc=1,\,\, a,b,c,d\in\BR\biggr\}$$ of $\rSL_2(\BR)^n$, which is an
embedded image of $\rSL_2(\BR)$ and denote it as $\rSL_2(\BR)$ with
abuse of notation. Thus the action of $\rSL_2(\BR)$ on $\BR^{2n}$
and a subgroup $\rSL_2(\BR)\ltimes\BR^{2n}$ of $\bG^n$ are
well-defined.

For $\bba=(a_1,\cdots,a_n)\in\BN^n$,
denote $$\rA_{\bba}=\biggl(\begin{pmatrix} \rD_{\sqrt\bba} & \rD_{\bold 0} \\
\rD_{\bold 0} & \rD_{\f{1}{\sqrt\bba}}
\end{pmatrix};\begin{pmatrix} \bold 0 \\
\bold 0
\end{pmatrix}\biggr)\in\bG^n$$ where $\sqrt\bba=(\sqrt
{a_1},\cdots,\sqrt{a_n})$,
$\f{1}{\sqrt\bba}=(\f{1}{\sqrt{a_1}},\cdots,\f{1}{\sqrt{a_n}})$. Let
$\rL_{\bold 1}=\rSL_2(\BR)\ltimes\BR^{2n}$ and
$$\rL_{\bba}=\biggl\{\biggl(\begin{pmatrix} a\,\rI_n & b\,\rD_{\bba} \\
c\,\rD_{\f{1}{\bba}} & d\,\rI_n
\end{pmatrix};\begin{pmatrix} \bbx \\
\bby
\end{pmatrix}\biggr):ad-bc=1,\,\,a,b,c,d\in\BR,\,\bbx,\bby\in\BR^n\biggr\}.$$

First, we assume that $a^*:=\gcd(a_1,\cdots,a_n)=1$. Then we
consider a lattice $\Gm^0_{\bba}$ of $\rL_{\bba}$ defined by
\begin{equation*}\begin{split}
\Gm^0_{\bba}&:=\biggl\{\biggl(\begin{pmatrix} a\,\rI_n & b\,\rD_{\bba} \\
c\,\rD_{\f{1}{\bba}} & d\,\rI_n
\end{pmatrix};\begin{pmatrix} \bbm \\
\bbh
\end{pmatrix}\biggr):ad-bc=1,\,\,a,b,d,\f{c}{a_*}\in\BZ,\,\bbm,\bbh\in\BZ^n\biggr\}
\end{split}\end{equation*}
where $a_*=\lcm(a_1,\cdots,a_n)$. Let $\Gm$ be any finite-indexed
subgroup of $\rSL_2(\BZ)^n\ltimes\BZ^{2n}$. Then we see that
$\Gm_{\bba}:=\Gm\cap\rL_{\bba}$ is a finite-indexed subgroup of
$\Gm^0_{\bba}$, and thus $\Gm_{\bba}$ is a lattice of $\rL_{\bba}$.
For each $\bba\in\BN^n$, define the {\it unipotent flow}
$\fU_{\bba}^u:\Gm\s\bG^n\to\Gm\s\bG^n$ by right translation by
$$\cU^u_{\bba}=\biggl(\begin{pmatrix} \rI_n & u\,\rD_{\bba} \\
\rD_{\bold 0} & \rI_n
\end{pmatrix};\begin{pmatrix} \bold 0 \\
\bold 0
\end{pmatrix}\biggr)\in\rL_{\bba}\subset\bG^n,$$ that is, $\fU_{\bba}^u(\Gm\gm)=\Gm\gm\cU^u_{\bba}$ for $u\in\BR$,
and the {\it geodesic flow} $\fG_{\bba}^t:\Gm\s\bG^n\to\Gm\s\bG^n$
by right translation by
$$\cG^t_0=\biggl(\begin{pmatrix} e^{-\f{t}{2}}\rI_n & \rD_{\bold 0} \\
\rD_{\bold 0} & e^{\f{t}{2}}\rI_n
\end{pmatrix};\begin{pmatrix} \bold 0 \\
\bold 0
\end{pmatrix}\biggr)\in\rL_{\bba}\subset\bG^n,$$ that is, $\fG_{\bba}^t(\Gm\gm)=\Gm\gm\cG^t_0$ for $t\in\BR$.
Note that if $a^*:=\gcd(a_1,\cdots,a_n)=1$, then $\cU^u_{\bba}$ and
$\cG^t_0$ satisfy the relation $\cU^u_{\bba}=\rA_{\bba}\cU^u_{\bold
1}\rA^{-1}_{\bba}$ and $\cG^t_0=\rA_{\bba}\cG^t_0\rA^{-1}_{\bba},$
where $\bold{1}=(1,\cdots,1)\in\BZ_+^n$ is the identity element of
the multiplicative semigroup $\BZ_+^n$. Since $\cU^u_{\bba}$ and
$\cG^t_0$ are in $\rL_{\bba}$, we note that {\it $\fU^u_{\bba}$ and
$\fG^t_{\bba}$ are the unipotent flow and geodesic flow on
$\Gm_{\bba}\s\rL_{\bba}$, respectively, as well}.

If $a^*\neq 1$, then we may write $\bba=a^*\bbb$ where
$\gcd(b_1,\cdots,b_n)=1$ and $a_k=a^* b_k$ for $k=1,\cdots,n$, and
also we observe that
$$\cU^{u/a^*}_{\bba}=\rA_{\bba}\cU^{u/a^*}_{\bold
1}\rA^{-1}_{\bba}=\rA_{\bbb}\cU^u_{\bold
1}\rA^{-1}_{\bbb}=\cU^u_{\bbb}$$ and
$$\cG^t_0=\rA_{\bba}\cG^t_0\rA^{-1}_{\bba}=\rA_{\bbb}\cG^t_0\rA^{-1}_{\bbb}.$$
In this occassion, we can use $\bbb$ in place of $\bba$. Thus,
without loss of generality, {\it we may assume that
$a^*:=\gcd(a_1,\cdots,a_n)=1$ in what follows}.

Let us consider for a moment that $\Gm_{\bba}=\Gm_{\bba}^0$.
Then for $$\gm_0=\biggl(\begin{pmatrix} \rI_n & \rD_{\bold 0} \\
\rD_{\bold 0} & \rI_n
\end{pmatrix};\begin{pmatrix} \bbx \\
\bold 0
\end{pmatrix}\biggr)\in\rL_{\bba}\subset\bG^n,\,\,\bbx\in\BR^n,$$ we have the relation
$\Gm_{\bba}^0\gm_0\cU^{u+1}_{\bba}\cG_0^t=\Gm_{\bba}^0\gm_0\cU^u_{\bba}\cG_0^t$
because $\cU^1_{\bba}\in\Gm_{\bba}^0$ and $\gm_0$ commutes with
$\cU^u_{\bba}$ by the multiplication law on $\bG^n$. Thus the closed
orbit $\{\Gm_{\bba}^0\gm_0\cU^u_{\bba}\cG_0^t:u\in\BR\}$ can be
described as $\{\Gm_{\bba}^0\gm_0\cU^u_{\bba}\cG_0^t:u\in[0,1)\}$
for each $t\in\BR$.

If $\Gm_{\bba}$ is a subgroup of $\Gm^0_{\bba}$ with finite index
$[\Gm^0_{\bba}:\Gm_{\bba}]$, then the manifold
$\Gm_{\bba}\s\rL_{\bba}$ is a finite covering of
$\Gm^0_{\bba}\s\rL_{\bba}$. Hence we can find some integer
$r=r(\Gm,\bba)\in [1,[\Gm^0_{\bba}:\Gm_{\bba}]]$ such that
$\{\Gm_{\bba}\gm_0\cU^u_{\bba}\cG_0^t:u\in[0,r)\}$ represents a
closed orbit $\{\Gm_{\bba}\gm_0\cU^u_{\bba}\cG_0^t:u\in\BR\}$ on
$\Gm_{\bba}\s\rL_{\bba}$.

Let $\Gm$ be any finite-indexed subgroup of
$\rSL_2(\BZ)^n\ltimes\BZ^{2n}$ and let $F$ be a $\Gm$-left invariant
function defined on $\bG^n$. Then its restriction $F|_{\rL_{\bba}}$
to $\rL_{\bba}$ is a $\Gm_{\bba}$-left invariant function defined on
$\rL_{\bba}$. We denote by $\overline{F|_{\rL_{\bba}}}$ the induced
function of $F|_{\rL_{\bba}}$ on $\Gm_{\bba}\s\rL_{\bba}$.

Then we prove the following theorem using Theorem 3.1 \cite{M3}.

\begin{thm} Let $F$ be a bounded continuous function on $\bG^n$
that is left invariant under a finite-indexed subgroup $\Gm$ of
$\rSL_2(\BZ)^n\ltimes\BZ^{2n}$, and let $$\gm_0=\biggl(\begin{pmatrix} \rI_n & \rD_{\bold 0} \\
\rD_{\bold 0} & \rI_n
\end{pmatrix};\begin{pmatrix} \bbx \\
\bold 0
\end{pmatrix}\biggr)\in\rL_{\bba}\subset\bG^n$$ where the components of the vector
$(\bbx,1)\in\BR^{n+1}$ are linearly independent over $\BQ$. If $h$
is a piecewise continuous function on $\BR/r\BZ$, then we have that
$$\lim_{t\to\iy}\f{1}{r}\int_0^r\overline
{F|_{\rL_{\bba}}}\circ\fG_{\bba}^t\circ\fU_{\bba}^u(\Gm_{\bba}\gm_0)h(u)\,du
=\int_{\Gm_{\bba}\s\rL_{\bba}}\overline{F|_{\rL_{\bba}}}
\,d\mu_{\bba}\,\f{1}{r}\int_0^rh(u)\,du$$ where $\mu_{\bba}$ is the
normalized Haar measure on $\Gm_{\bba}\s\rL_{\bba}$ and
$\Gm_{\bba}=\Gm\cap\rL_{\bba}$.
\end{thm}

\pf We consider a map $\Xi$ from $\rL_{\bba}$ into $\rL_{\bold
1}=\rSL_2(\BR)\ltimes\BR^{2n}$ defined by
$$\Xi\biggl[\biggl(\begin{pmatrix} a\,\rI_n & b\,\rD_{\bba} \\
c\,\rD_{\f{1}{\bba}} & d\,\rI_n
\end{pmatrix};\begin{pmatrix} \bbx \\
\bby
\end{pmatrix}\biggr)\biggr]=\biggl(\begin{pmatrix} a\,\rI_n & b\,\rI_n \\
c\,\rI_n & d\,\rI_n
\end{pmatrix};\begin{pmatrix} \bbx \\
\bba\bby
\end{pmatrix}\biggr)$$ where $\bba\bby=(a_1 y_1,\cdots,a_n y_n)$.
Then it is easy to check that $\Xi$ is bijective and a group
homomorphism. Thus $\,\Xi:\rL_{\bba}\to\rL_{\bold 1}$ is a group
isomorphism, $F|_{\rL_{\bba}}\circ\Xi^{-1}$ is a bounded continuous
function on $\rL_{\bold 1}$ which is left invariant under
$\Xi(\Gm_{\bba})$, and the induced map
$\overline\Xi:\Gm_{\bba}\s\rL_{\bba}\to\Xi(\Gm_{\bba})\s\rL_{\bold
1}$ defined by $\overline\Xi(\Gm_{\bba}\gm)=\Xi(\Gm_{\bba})\Xi(\gm)$
is a group isomorphism. We observe that
$$\overline{{F|_{\rL_{\bba}}\circ\Xi^{-1}}}=\overline{F|_{\rL_{\bba}}}\circ{\overline\Xi}^{\,-1}.$$

Since $\Xi(\Gm_{\bba}^0)$ must be a finite-indexed subgroup of
$\rSL_2(\BZ)\ltimes\BZ^{2n}$, we see that $\Xi(\Gm_{\bba})$ is also
a finite-indexed subgroup of $\rSL_2(\BZ)\ltimes\BZ^{2n}$. Hence it
follows from Theorem 3.1 \cite{M3} that
\begin{equation*}\begin{split}&\lim_{t\to\iy}\f{1}{r}\int_0^r\overline
{(F|_{\rL_{\bba}}\circ\Xi^{-1})}\circ\fG_{\bold 1}^t\circ\fU_{\bold
1}^u(\Xi(\Gm_{\bba})\gm_0)h(u)\,du\\&\qquad\qquad\qquad\qquad=\int_{\Xi(\Gm_{\bba})\s\rL_{\bold
1}}\overline{F|_{\rL_{\bba}}\circ\Xi^{-1}} \,d\mu_{\bold
1}\,\f{1}{r}\int_0^rh(u)\,du
\end{split}\end{equation*} where $\mu_{\bold 1}$ is the
normalized Haar measure on $\Xi(\Gm_{\bba})\s\rL_{\bold 1}$. We note
that $\Xi(\cU^u_{\bba})=\cU^u_{\bold 1}$, $\Xi(\cG_0^t)=\cG_0^t$ and
$\Xi(\gm_0)=\gm_0$. Since $\Xi$ is a group isomorphism, by using the
change of variables we can conclude that
$$\lim_{t\to\iy}\f{1}{r}\int_0^r\overline
{F|_{\rL_{\bba}}}\circ\fG_{\bba}^t\circ\fU_{\bba}^u(\Gm_{\bba}\gm_0)h(u)\,du
=\int_{\Gm_{\bba}\s\rL_{\bba}}\overline{F|_{\rL_{\bba}}}
\,d\mu_{\bba}\,\f{1}{r}\int_0^rh(u)\,du$$ where $\mu_{\bba}$ is the
normalized Haar measure on $\Gm_{\bba}\s\rL_{\bba}$. Therefore we
are done. \qed

Using the isomorphism $\Xi$ between $\rL_{\bba}$ and $\rL_{\bold
1}=\rSL_2(\BR)\ltimes\BR^{2n}$, we may extend the above theorem to
dominated unbounded functions. Let $F$ be a function on $\bG^n$
which is left invariant under a finite-indexed subgroup $\Gm$ of
$\rSL_2(\BZ)^n\ltimes\BZ^{2n}$. Then we say that {\it $F$ is
dominated by $F_R$ on $\rL_{\bba}$}, if there is a constant $L>1$
such that
$$|\,F|_{\rL_{\bba}}\circ\Xi^{-1}(z\bold 1,\phi\bold 1;\bxi)\,|X_R(z)\le
L+F_R(z;\bxi)$$ for any sufficiently large $R>1$, uniformly for all
$(z\bold 1,\phi\bold 1;\bxi)\in\rL_{\bold 1}$. Here as in
\cite{M1,M2,M3} we see that
$$X_R(z)=\sum_{\gm\in\pm\Gm_{\iy}\s\rSL_2(\BZ)}\mathbbm{1}_{[R,\iy)}(v_{\gm})$$
and
$$F_R(z;\bxi)=\sum_{\gm\in\Gm_{\iy}\s\rSL_2(\BZ)}\sum_{\bbm\in\BZ^n}f\bigl((\bby_{\gm}+\bbm)v_{\gm}^{\f{1}{2}}\bigr)
v_{\gm}^{\f{n}{2}}\mathbbm{1}_{[R,\iy)}(v_{\gm})$$ for $f\in
C(\BR^n)$ with rapid decay, where $v_{\gm}>0$ and
$\bby_{\gm}\in\BR^n$ are defined by $$v_{\gm}:=\Ima(\gm
z)=\Ima\biggl(\f{az+b}{cz+d}\biggr)=\f{v}{|cz+d|^2}\,\,\text{ for
$z=u+i v$, }$$
$$\begin{pmatrix} \bbx_{\gm} \\
\bby_{\gm}
\end{pmatrix}=\gm\begin{pmatrix} \bbx \\
\bby
\end{pmatrix}\text{ for $\gm=\begin{pmatrix} a & b \\
c & d
\end{pmatrix}$,   and }
\Gm_{\iy}=\biggl\{\begin{pmatrix} \rI_n & m\rI_n \\
\rD_{\bold 0} & \rI_n
\end{pmatrix}:m\in\BZ\biggr\}.$$

We observe that if $F$ is a $\Gm$-left invariant function on $\bG^n$
where $\Gm$ is a finite-indexed subgroup of
$\rSL_2(\BZ)^n\ltimes\BZ^{2n}$, then the function $F_{\bba}$ defined
by
\begin{equation}F_{\bba}(\gm)=F(\gm\rA_{\bba})\,\,\text{ for $\gm\in\bG^n$ and
$\bba\in\BR^n$ }
\end{equation} is $\Gm$-left invariant on $\bG^n$ and $F_{\bba}|_{\rL_{\bba}}$ is $\Gm_{\bba}$-left invariant on $\rL_{\bba}$.

\begin{thm} Let $F\ge 0$ be a continuous function on $\bG^n$
which is left invariant under a finite-indexed subgroup $\Gm$ of
$\rSL_2(\BZ)^n\ltimes\BZ^{2n}$, and let $$\gm_0=\biggl(\begin{pmatrix} \rI_n & \rD_{\bold 0} \\
\rD_{\bold 0} & \rI_n
\end{pmatrix};\begin{pmatrix} \bbx \\
\bold 0
\end{pmatrix}\biggr)\in\rL_{\bba}\subset\bG^n,\,\bba\in\BN^n,$$ where $\bbx\in\BT^n$ is a vector of diophantine type
$\kp<\f{d-1}{d-2}$
and the components of the vector $(\bbx^t,1)\in\BR^{n+1}$ are
linearly independent over $\BQ$. If $\,h$ is a piecewise continuous
function on $\BR/r\BZ$ and $F$ is dominated by $F_R$ on
$\rL_{\bba}$, then we have that
$$\lim_{v\to 0}\f{1}{r}\int_0^r \overline{{F_{\bba}|_{\rL_{\bba}}}\circ\Xi^{-1}}\biggl(u+i v,\bold 0;\begin{pmatrix} \bbx \\
\bold 0
\end{pmatrix}\biggr)h(u)\,du=\int_{\Gm_{\bba}\s\rL_{\bba}}\overline{F_{\bba}|_{\rL_{\bba}}}
\,d\mu_{\bba}\,\f{1}{r}\int_0^rh(u)\,du$$ where $\mu_{\bba}$ is the
normalized Haar measure on $\Gm_{\bba}\s\rL_{\bba}$.
\end{thm}

\pf Since we see that
\begin{equation*}\begin{split}\int_{\Gm_{\bba}\s\rL_{\bba}}\overline{F_{\bba}|_{\rL_{\bba}}}
\,d\mu_{\bba}&=\int_{\Xi(\Gm_{\bba})\s\rL_{\bold
1}}\overline{F_{\bba}|_{\rL_{\bba}}\circ\Xi^{-1}} \,d\mu_{\bold
1}\\&=\int_{\Xi(\Gm_{\bba})\s\rL_{\bold
1}}\overline{F_{\bba}|_{\rL_{\bba}}}\circ\overline{\Xi}^{\,-1}
\,d\mu_{\bold 1}\end{split}\end{equation*} where $\mu_{\bold 1}$ is
the normalized Haar measure on $\Xi(\Gm_{\bba})\s\rL_{\bold 1}$, it
follows from Theorem 5.1 in \cite{M3} and its little modification.
\qed

We recall that
$$\Gm^n=\biggl\{\biggl(\rM;\f{1}{2}\begin{pmatrix} \bba\bbb \\
\bbc\bbd
\end{pmatrix}+\begin{pmatrix} \bbm \\
\bbh
\end{pmatrix}\biggr):\bba,\bbb,\bbc,\bbd,\bbm,\bbh\in\BZ^n,\rM=\begin{pmatrix} \rD_{\bba} & \rD_{\bbb} \\
\rD_{\bbc} & \rD_{\bbd} \end{pmatrix}\in\rSL_2(\BZ)^n\biggr\}
$$ where $\bba\bbb=(a_1 b_1,\cdots,a_n b_n)$ and $\bbc\bbd=(c_1
d_1,\cdots,c_n d_n)$. Then we see that $\Gm^n$ is of finite index in
$\rSL_2(\BZ)^n\ltimes(\f{1}{2}\BZ)^{2n}$. So we remark that we can
extend the above theorems to
$\rSL_2(\BZ)^n\ltimes(\f{1}{2}\BZ)^{2n}$ rather than
$\rSL_2(\BZ)^n\ltimes\BZ^{2n}$. We set
\begin{equation*}\begin{split}
\Gm^*_{\bba}&:=\biggl\{\biggl(\begin{pmatrix} a\,\rI_n & b\,\rD_{\bba} \\
c\,\rD_{\f{1}{\bba}} & d\,\rI_n
\end{pmatrix};\f{1}{2}\begin{pmatrix} \bbm \\
\bbh\end{pmatrix}\biggr):ad-bc=1,\,\,a,b,d,\f{c}{a_*}\in\BZ,\,\bbm,\bbh\in\BZ^n\biggr\}
\end{split}\end{equation*} and
\begin{equation*}
\Gm^n_{\bba}=\biggl\{\biggl(\begin{pmatrix} a\,\rI_n & b\,\rD_{\bba} \\
c\,\rD_{\f{1}{\bba}} & d\,\rI_n
\end{pmatrix};\f{1}{2}\begin{pmatrix} ab\,\bba \\
cd\f{1}{\bba}
\end{pmatrix}+\begin{pmatrix} \bbm \\
\bbh
\end{pmatrix}\biggr):ad-bc=1,
a,b,d,\f{c}{a_*}\in\BZ,\bbm,\bbh\in\BZ^n\biggr\}
\end{equation*} where $a_*=\lcm(a_1,\cdots,a_n)$.
Then we see that $\Gm^*_{\bba}$ is a lattice of $\rL_{\bba}$, and
also we use $\Gm^*_{\bba}$ instead of $\Gm^0_{\bba}$ in the
extension of the above theorems to
$\rSL_2(\BZ)^n\ltimes(\f{1}{2}\BZ)^{2n}$. It is easy to check that
$\Gm_{\bba}^n$ is a finite-indexed subgroup of $\Gm^*_{\bba}$, and
so it is a lattice of $\rL_{\bba}$. Thus we see that
$\Xi(\Gm^n_{\bba})\subset\rSL_2(\BZ)\ltimes(\f{1}{2}\BZ)^{2n}\subset\rL_{\bold
1}$, $\Xi(\Gm^n_{\bba})$ is a lattice of $\rL_{\bold 1}$ and
\begin{equation*}
\Xi(\Gm^n_{\bba})=\biggl\{\biggl(\begin{pmatrix} a\,\rI_n & b\,\rI_n \\
c\,\rI_n & d\,\rI_n
\end{pmatrix};\f{1}{2}\begin{pmatrix}ab\,\bba \\
cd\,\bold 1
\end{pmatrix}+\begin{pmatrix} \bbm \\
\bba\bbh
\end{pmatrix}\biggr):\begin{pmatrix} a\,\rI_n & b\,\rI_n \\
c\,\rI_n & d\,\rI_n
\end{pmatrix}\in\Gm_0(a_*),\bbm,\bbh\in\BZ^n\biggr\}.
\end{equation*} where $\Gm_0(a_*)$ is a congruence subgroup of
$\rSL_2(\BZ)$ given by
\begin{equation*}
\Gm_0(a_*):=\biggl\{\biggl(\begin{pmatrix} a\,\rI_n & b\,\rI_n \\
c\,\rI_n & d\,\rI_n
\end{pmatrix}:ad-bc=1,
a,b,d,\f{c}{a_*}\in\BZ\biggr\}.
\end{equation*}

\begin{prop} The left action of the group $\Xi(\Gm^n_{\bba})$ is properly
discontinuous. In addition, a fundamental domain
$\Pi_{\Xi(\Gm^n_{\bba})}$ of $\Xi(\Gm^n_{\bba})$ in $\rL_{\bold
1}=\rSL_2(\BR)\ltimes\BR^{2n}$ is given by
$$\cF_{\Gm_0(a_*)}\times[0,\pi)\times\bigl[-\f{1}{2},\f{1}{2}\bigr)^n
\times\prod_{k=1}^n\bigl[-\f{1}{2}a_k,\f{1}{2}a_k\bigr)$$ where
$\cF_{\Gm_0(a_*)}$ is the fundamental domain in $\BH^2$ of the group
$\Gm_0(a_*)$.
\end{prop}

\pf We take any matrix $\biggl(\begin{pmatrix} p\,\rI_n & q\,\rI_n \\
r\,\rI_n & s\,\rI_n
\end{pmatrix};\begin{pmatrix} \bzt \\
\bta
\end{pmatrix}\biggr)\in\rL_{\bold 1}$. Then there exists a unique
matrix
$\begin{pmatrix} a\,\rI_n & b\,\rI_n \\
c\,\rI_n & d\,\rI_n
\end{pmatrix}\in\cF_{\Gm_0(a_*)}$ such that
$$\begin{pmatrix} a\,\rI_n & b\,\rI_n \\
c\,\rI_n & d\,\rI_n
\end{pmatrix}\begin{pmatrix} p\,\rI_n & q\,\rI_n \\
r\,\rI_n & s\,\rI_n
\end{pmatrix}\in\cF_{\Gm_0(a_*)}.$$ Thus we may choose some vectors
$\bbm,\bbh\in\BZ^n$ so that the element
\begin{equation*}\begin{split}&\biggl(\begin{pmatrix} a\,\rI_n & b\,\rI_n \\
c\,\rI_n & d\,\rI_n
\end{pmatrix};\f{1}{2}\begin{pmatrix} ab\bba \\
cd\bold 1
\end{pmatrix}+\begin{pmatrix} \bbm \\
\bba\bbh
\end{pmatrix}\biggr)\biggl(\begin{pmatrix} p\,\rI_n & q\,\rI_n \\
r\,\rI_n & s\,\rI_n
\end{pmatrix};\begin{pmatrix} \bzt \\
\bta
\end{pmatrix}\biggr)\\
&=\biggl(\begin{pmatrix} a\,\rI_n & b\,\rI_n \\
c\,\rI_n & d\,\rI_n
\end{pmatrix}\begin{pmatrix} p\,\rI_n & q\,\rI_n \\
r\,\rI_n & s\,\rI_n
\end{pmatrix};\f{1}{2}\begin{pmatrix} ab\bba \\
cd\bold 1
\end{pmatrix}+\begin{pmatrix} \bbm \\
\bba\bbh
\end{pmatrix}+\begin{pmatrix} a\,\rI_n & b\,\rI_n \\
c\,\rI_n & d\,\rI_n
\end{pmatrix}\begin{pmatrix} \bzt \\
\bta
\end{pmatrix}\biggr)
\end{split}\end{equation*} is contained in
$$\cF_{\Gm_0(a_*)}\times[0,\pi)\times\bigl[-\f{1}{2},\f{1}{2}\bigr)^n
\times\prod_{k=1}^n\bigl[-\f{1}{2}a_k,\f{1}{2}a_k\bigr).$$ Hence we
complete the proof. \qed

We now consider
\begin{equation}F(\bbz,\bphi;\bxi)=\Th_f(\bbz,\bphi;\bxi)\overline{\Th_g(\bbz,\bphi;\bxi)}\end{equation}
for $f,g\in\cS(\BR^n)$, where
$\Th_f(\bbz,\bphi;\bxi)=\Th_f(\bbz,\bphi;\bxi,0)$. Using the group
isomorphism $\Xi$, we can regard
$F_{\bba}|_{\rL_{\bba}}\circ\Xi^{-1}$ as a function on $\rL_{\bold
1}=\rSL_2(\BR)\ltimes\BR^{2n}$. For $\bba=(a_1,\cdots,a_n)\in\BN^n$,
recall the notation denoted by
$$\rA_{\bba}=\biggl(\rM_{\bba};\begin{pmatrix} \bold 0 \\
\bold 0
\end{pmatrix}\biggr)\in\bG^n,\,\,\,\rM_{\bba}=\begin{pmatrix} \rD_{\sqrt\bba} & \rD_{\bold 0} \\
\rD_{\bold 0} & \rD_{\f{1}{\sqrt\bba}}
\end{pmatrix}$$ where $\sqrt\bba=(\sqrt
{a_1},\cdots,\sqrt{a_n})$,
$\f{1}{\sqrt\bba}=(\f{1}{\sqrt{a_1}},\cdots,\f{1}{\sqrt{a_n}})$.
From the above observation, one sees that the map
$\Xi:\rL_{\bba}\to\rL_{\bold 1}$ given by
\begin{equation}\Xi(\bbz,\bphi;\bxi)=\biggl(\rM_{\bba}^{-1}(\bbz,\bphi)\rM_{\bba};\begin{pmatrix} \bbx \\
\bba\bby
\end{pmatrix}\biggr),\,\,\bxi=\begin{pmatrix} \bbx \\
\bby
\end{pmatrix},\end{equation} is a group isomorphism. By the unique Iwasawa
decomposition (3.6), we note that
\begin{equation}\begin{split}&\rM_{\bba}(z\bold 1,\phi\bold 1)\\&=\rM_{\bba}\begin{pmatrix} \rI_n & u\,\rI_n \\
\rD_{\bold 0} & \rI_n \end{pmatrix}\rM_{\bba}^{-1}\rM_{\bba}\begin{pmatrix} v^{\f{1}{2}}\rI_n  & \rD_{\bold 0} \\
\rD_{\bold 0} & v^{-\f{1}{2}}\rI_n \end{pmatrix}\begin{pmatrix} \cos\phi\,\rI_n & -\sin\phi\,\rI_n \\
\sin\phi\,\rI_n & \cos\phi\,\rI_n
\end{pmatrix}\\
&=\begin{pmatrix} \rI_n & u\rD_{\bba} \\
\rD_{\bold 0} & \rI_n \end{pmatrix}\begin{pmatrix} v^{\f{1}{2}}\rD_{\sqrt\bba}  & \rD_{\bold 0} \\
\rD_{\bold 0} & v^{-\f{1}{2}}\rD_{1/\sqrt\bba} \end{pmatrix}\begin{pmatrix} \cos\phi\,\rI_n & -\sin\phi\,\rI_n \\
\sin\phi\,\rI_n & \cos\phi\,\rI_n\end{pmatrix}
\end{split}\end{equation} for any $(z\bold 1,\phi\bold
1)\in\rSL_2(\BR)$. Thus by the definition of $\Th_f$ and (3.18), for
$(z\bold 1,\phi\bold 1;\bxi)\in\rSL_2(\BR)\ltimes\BR^{2n}$, we have
that
\begin{equation}\begin{split}&F_{\bba}|_{\rL_{\bba}}\circ\Xi^{-1}(z\bold 1,\phi\bold
1;\bxi)\\&=\Th_f\biggl(\biggl(\rM_{\bba}(z\bold
1,\phi\bold 1)\rM_{\bba}^{-1};\begin{pmatrix} \bbx \\
\f{1}{\bba}\bby
\end{pmatrix}\biggr)\rA_{\bba}\biggr)\overline{\Th_g\biggl(\biggl(\rM_{\bba}(z\bold
1,\phi\bold 1)\rM_{\bba}^{-1};\begin{pmatrix} \bbx \\
\f{1}{\bba}\bby
\end{pmatrix}\biggr)\rA_{\bba}\biggr)}\\
&=\Th_f\biggl(\rM_{\bba}(z\bold
1,\phi\bold 1);\begin{pmatrix} \bbx \\
\f{1}{\bba}\bby
\end{pmatrix}\biggr)\overline{\Th_g\biggl(\rM_{\bba}(z\bold
1,\phi\bold 1);\begin{pmatrix} \bbx \\
\f{1}{\bba}\bby
\end{pmatrix}\biggr)}\\
&=(a_1\cdots
a_n)^{\f{1}{2}}v^{\f{n}{2}}\sum_{\bbm,\bbh\in\BZ^n}f_{\phi\bold
1}\biggl((m_1-\f{1}{a_1}y_1)\sqrt{a_1
v},\cdots,(m_n-\f{1}{a_n}y_n)\sqrt{a_n v}\biggr)\\
&\qquad\qquad\qquad\qquad\quad\times\overline{g_{\phi\bold
1}\biggl((h_1-\f{1}{a_1}y_1)\sqrt{a_1
v},\cdots,(h_n-\f{1}{a_n}y_n)\sqrt{a_n v}\biggr)}\\
&\qquad\qquad\qquad\quad\times\fe(\bbh\cdot\bbx-\bbm\cdot\bbx)
\,\fe\biggl(\f{1}{2}u\sum_{k=1}^n a_k
\bigl((m_k-\f{1}{a_k}y_k)^2-(h_k-\f{1}{a_k}y_k)^2\bigr)\biggr)\\
&=\Th_f(\bba z,\phi\bold 1;\bxi)\,\overline{\Th_g(\bba z,\phi\bold
1;\bxi)}
\end{split}\end{equation}

\begin{prop} If $f,g\in\cS(\BR^n)$ and $\bba\in\BN^n$, then $F_{\bba}$ is dominated by
$F_R$ on $\rL_{\bba}$.
\end{prop}

\pf We now set $$f^*(\bbw)=\sqrt{a_1\cdots
a_n}\,\,\sup_{\phi\in\BR}\biggl|f_{\phi\bold
1}\biggl(\f{1}{\sqrt\bba}\,\bbw\biggr)g_{\phi\bold
1}\biggl(\f{1}{\sqrt\bba}\,\bbw\biggr)\biggr|.$$ Since
$f_{\bphi}\in\cS(\BR^n)$ for $f\in\cS(\BR^n)$, we see that for any
$R>1$ there is a constant $C_R>0$ such that
\begin{equation}f^*(\bbw)\le\sqrt{a_1\cdots
a_n}\,\,C_R\left(1+|\bbw|_{1/\bba}\,\right)^{-2R}
\end{equation} for any $\bbw\in\BR^n$. So the function $f^*$ has
rapid decay. Let us define
$$F_R(z;\bxi)=\sum_{\gm\in\Gm_{\iy}\s\rSL_2(\BZ)}\sum_{\bbh\in\BZ^n}f\bigl((\bby_{\gm}+\bbh)v_{\gm}^{\f{1}{2}}\bigr)
v_{\gm}^{\f{n}{2}}\mathbbm{1}_{[R,\iy)}(v_{\gm}).$$ If $v>R$ for
sufficiently large $R>1$, then as in \cite{M1,M2,M3} we have that
$$F_R(z;\bxi)=v^{\f{n}{2}}\sum_{\bbh\in\BZ^n}
\biggl[f^*\biggl((\bbh+\bby)v^{\f{1}{2}}\biggr)+f^*\biggl((\bbh-\bby)v^{\f{1}{2}}\biggr)\biggr].$$
From (3.14) and \cite{M1,M2}, we note that
\begin{equation*}\begin{split}&F_{\bba}|_{\rL_{\bba}}\circ\Xi^{-1}(z\bold 1,\phi\bold
1;\bxi)\\&\quad=\sqrt{a_1\cdots a_n}\,v^{\f{n}{2}}f_{\phi\bold
1}\biggl(\bigl(\bbh-\f{1}{\bba}\,\bby\bigr)\sqrt\bba\,\,v^{\f{1}{2}}\biggr)\overline{g_{\phi\bold
1}\biggl(\bigl(\bbh-\f{1}{\bba}\,\bby\bigr)\sqrt\bba\,\,v^{\f{1}{2}}\biggr)}+\cO(v^{-R})
\end{split}\end{equation*}
uniformly for all $\bby\in\BR^n$ such that
$h_k-\f{1}{2}\le\f{1}{a_k}\,y_k< h_k+\f{1}{2}$ for all
$k=1,\cdots,n$. Thus by the definition of $f^*$, we obtain that
\begin{equation*}\begin{split}&|\,F_{\bba}|_{\rL_{\bba}}\circ\Xi^{-1}(z\bold 1,\phi\bold
1;\bxi)\,|\\&=\sqrt{a_1\cdots a_n}\,v^{\f{n}{2}}\biggl|f_{\phi\bold
1}\biggl(\bigl(\bbh-\f{1}{\bba}\,\bby\bigr)\sqrt\bba\,\,v^{\f{1}{2}}\biggr)\overline{g_{\phi\bold
1}\biggl(\bigl(\bbh-\f{1}{\bba}\,\bby\bigr)\sqrt\bba\,\,v^{\f{1}{2}}\biggr)}\,\biggr|+\cO(v^{-R})\\
&\le f^*\bigl((\bba\bbh-\bby)v^{\f{1}{2}}\bigr)+\cO(v^{-R})\le
F_R(z;\bxi)+\cO(v^{-R}),
\end{split}\end{equation*} because $\bba\bbh\in\BZ^n$. Hence there is
a constant $L>1$ such that
$$|\,F_{\bba}|_{\rL_{\bba}}\circ\Xi^{-1}(z\bold 1,\phi\bold
1;\bxi)\,|X_R(z)\le L+F_R(z;\bxi)$$ for any sufficiently large
$R>1$, uniformly for all $(z\bold 1,\phi\bold 1;\bxi)\in\rL_{\bold
1}$. \qed

\begin{prop} For $\bba\in\BN^n$, let $F_{\bba}$ be the function
defined as in $(3.15)$ where $F$ is the continuous function on
$\bG^n$ as in $(3.16)$. Then we have that
$$\int_{\Gm_{\bba}\s\rL_{\bba}}\overline{F_{\bba}|_{\rL_{\bba}}}
\,d\mu_{\bba}=\int_{\BR^n}f(\bby)\overline{g(\bby)}\,d\bby$$ where
$\mu_{\bba}$ is the normalized Haar measure on
$\Gm_{\bba}\s\rL_{\bba}$.
\end{prop}

\pf Since $\Gm_{\bba}^n\subset\Gm^n$, by Proposition 3.6 we see that
$F$ is left invariant under a finite-indexed subgroup $\Gm_{\bba}^n$
of $\Gm_{\bba}^*$, and so is $F_{\bba}$. By applying the extension
of Theorem 3.8 to $\rSL_2(\BZ)^n\ltimes(\f{1}{2}\BZ)^{2n}$ with
$\Gm_{\bba}=\Gm_{\bba}^n$ and $\Gm_{\bba}^0=\Gm_{\bba}^*$, we can
obtain that
$$\int_{\Gm_{\bba}^n\s\rL_{\bba}}\overline{F_{\bba}|_{\rL_{\bba}}}
\,d\mu_{\bba}=\int_{\Xi(\Gm_{\bba}^n)\s\rL_{\bold
1}}\overline{F_{\bba}|_{\rL_{\bba}}\circ\Xi^{-1}} \,d\mu_{\bold 1}$$
where $\mu_{\bold 1}$ is the normalized Haar measure on
$\Xi(\Gm_{\bba}^n)\s\rL_{\bold 1}$. We observe that
$$d\mu_{\bold 1}(\bbz,\bphi;\bxi)=\f{1}{\ld_0\pi(a_1\cdots
a_n)}\,\f{du\,dv\,d\phi\,d\bbx\,d\bby}{v^2},$$ since we see that
$\ld_0:=\int_{\cF_{\Gm_0(a_*)}}v^{-2}\,du\,dv>0$ by Proposition
3.10. Since $\widetilde\rR(i\bold 1,\phi\bold 1)$ is a unitary
operator and we know the fact that
$$\int_{[-\f{1}{2},\f{1}{2})^n}\fe(\bbh\cdot\bbx-\bbm\cdot\bbx)\,d\bbx=0$$
for any $\bbm,\bbh\in\BZ^n$ with $\bbm\neq\bbh$, it follows from
(3.19) and simple calculation that
\begin{equation*}\begin{split}&\int_{\Xi(\Gm_{\bba}^n)\s\rL_{\bold
1}}\overline{F_{\bba}|_{\rL_{\bba}}\circ\Xi^{-1}} \,d\mu_{\bold
1}=\f{(a_1\cdots a_n)^{\f{1}{2}}}{\ld_0\pi(a_1\cdots a_n)}
\int_0^{\pi}\int_{\cF_{\Gm_0(a_*)}}v^{\f{n}{2}-2}\\
&\quad\int_{\prod_{k=1}^n[-\f{1}{2}a_k,\f{1}{2}a_k)}\sum_{\bbh\in\BZ^n}f_{\phi\bold
1}\biggl((m_1-\f{1}{a_1}y_1)\sqrt{a_1
v},\cdots,(m_n-\f{1}{a_n}y_n)\sqrt{a_n v}\biggr)\\
&\qquad\qquad\qquad\times\overline{g_{\phi\bold
1}\biggl((h_1-\f{1}{a_1}y_1)\sqrt{a_1
v},\cdots,(h_n-\f{1}{a_n}y_n)\sqrt{a_n
v}\biggr)}\,d\bby\,dv\,du\,d\phi\\
&\qquad\quad=\f{1}{\ld_0\pi}\int_0^{\pi}\int_{\cF_{\Gm_0(a_*)}}\f{1}{v^2}
\biggl(\int_{\BR^n}f_{\phi\bold 1}(\bby)\overline{g_{\phi\bold
1}(\bby)}\,d\bby\biggr)dv\,du\,d\phi\\
&\qquad\quad=\f{1}{\ld_0\pi}\biggl(\int_{\BR^n}f(\bby)\overline{g(\bby)}\,d\bby\biggr)
\biggl(\int_0^{\pi}\int_{\cF_{\Gm_0(a_*)}}\f{1}{v^2}\,du\,dv\,d\phi\biggr)=\int_{\BR^n}f(\bby)\overline{g(\bby)}\,d\bby.\qed
\end{split}\end{equation*}

By Proposition 3.12, (3.19) and the extension of Theorem 3.9 to
$\rSL_2(\BZ)^n\ltimes(\f{1}{2}\BZ)^{2n}$ with
$\Gm_{\bba}=\Gm_{\bba}^n$ and $\Gm_{\bba}^0=\Gm_{\bba}^*$, we easily
obtain the following corollary.

\begin{cor}
Suppose that $f(\mathbf w) = \psi(|\mathbf w|^2)$ with $\psi \in \cS
(\R_+)$ real-valued, and let $h:\BR/2\BZ \to \BR_+$ be a nonnegative
piecewise continuous function. If $\bap\in\BT^n$ is a vector of
diophantine type $\kappa <\f{n-1}{n-2}$ so that the components of
the vector $(\bap,1)$ are linearly independent over $\BQ$, then we
have that
$$ \lim_{v \to 0} \frac{1}{2}\int_0^2 \biggl| \Theta_f \biggl(\bba z,\bold 0\,;\begin{pmatrix}
\bap \\ \bold 0 \end{pmatrix}\biggr)\biggr|^2 h(u)\,du =
\frac{n}{2}\,|B|\int_0^\infty \psi(r)^2 r^{n/2-1} dr
\frac{1}{2}\int_0^2 h(u)\,du,$$ where $z=u+i v\in\BH^2$ and
$\bba\in\BN^n$.
\end{cor}

{\bf Proof of Theorem 1.3.} We choose $f(\bbw)=g(\bbw)$ in Corollary
3.13 to be a rapidly decreasing function $\psi(|\mathbf w|^2)$ which
approximates the characteristic function
$\mathbbm{1}_{[0,1]}(|\bbw|^2)$. Take any $z=u+i v\in\BH^2$. Since
$f_0=f$, if we set $N=[1/v]$, then by (3.14) we have that
\begin{equation}\lim_{N\to\iy}\f{1}{N^{n/2}}\sum_{p=1}^N|\fr[\rD_{\bba},\bap](p)|^2=(a_1\cdots
a_n)^{-\f{1}{2}}\lim_{v\to 0}\f{1}{2}\int_0^2\bigl|\Theta_f(\bba
z,\bold 0\,; \bxi)\bigr|^2\,du.
\end{equation}
Thus by applying Corollary 3.13 with $h=1$, we conclude that
$$(a_1\cdots a_n)^{\f{1}{2}}\lim_{N \to
\infty}\frac{1}{N^{n/2}} \sum_{p=1}^{[N]}
|\fr[\rD_{\bba},\bap](p)|^2 = \frac{n}{2}\,|B| \int_0^\infty
\psi(r)^2 r^{n/2-1} dr,$$ where $|B|$ is the volume of the unit ball
in $\BR^n$. Hence we obtain our result because
$|E^{\rD_{\bba}}|=(a_1\cdots a_n)^{-1/2}|B|$. \qed

\section{ Proof of Theorem 1.2 and Theorem 1.4 }

The proof of Theorem 1.2 relies on the following proposition whose
proof can be found in \cite{KS}. For the proof of Theorem 1.4, we
require Theorem 4.4, Lemma 4.5 and Lemma 4.6 in which Property 1 is
assumed to hold.

\begin{prop} Let $\bap\in\BR^n$ $($$n\ge 2$$)$ and $\rM\in\fS_n^+(\BZ)$.
If $h:(0,\iy)\to\BR$ is a smooth function satisfying that
\begin{equation}\biggl|\biggl(\f{d}{dt}\biggr)^k
h(t)\biggr|\le\f{C_k}{t^k},\,\,t>1,
\end{equation} for all $k=0,1,\cdots\,,$ then for each
$\el=1,\cdots,N$, we have that
$$|\la F^{K,\el}_{\rM}(\,\cdot+s)h\ra_T|\ll T^{-N}$$ uniformly in $K\ge 1$ and $s\ge 0$, for any
$N>0$, where the implied constant depends only on $C_k$.
\end{prop}

At the final stage of the proof of Theorem 1.2 and Theorem 1.4, we
will show that $\rF^K_{\rM}(t)$ is a good approximation of
$\rF_{\rM}(t)$.

\begin{lemma} Let $\bap\in\BR^n$ $($$n\ge 2$$)$ and $\rM\in\fS_n^+(\BZ)$.
For any $T, K\ge 1$, we have that
$$\la|F_{\rM}(\,\cdot+s)-F_{\rM}^K(\,\cdot+s)|^2\ra_T\ll \f{T^{n-1}}{K}$$
uniformly in $s$ with $0\le s\ll T.$
\end{lemma}

\pf Since $\mu$ is supported in $(0,\iy)$, this result can be
achieved as in Lemma 2.8 \cite{KS} by replacing the Euclidean norm
$|\cdot|$ by the elliptic norm $|\cdot|_{\rM}$ as follows;
\begin{equation}|F_{\rM}(t+s)-F_{\rM}^K(t+s)|\le
t^{-\f{n-1}{2}}\sum_{\bbm\in\cA}1,\end{equation} where
$\cA=\{\bbm\in\BZ^n:|\,\,|\bbm-\bap|_{\rM}-(t+s)|\le K^{-1}\}$.
Hence this leads to the required estimate. \qed

{\bf Proof of Theorem 1.2.} As in Theorem 1.1 \cite{KS}, it easily
follows from (2.8), (2.9), Proposition 4.1 and Lemma 4.2.\qed

We furnish a useful lemma whose proof follows from the Abel's
summation formula as in Lemma 2.4 \cite{KS}.

\begin{lemma} Let $\rM\in\fS^+_n(\BZ)$ and $\bap\in\BR^n$ be given.
If there is some $a>0$ such that $\rR[\rM,\bap](p)\ll p^a$ for all
$p\in\BN$, then for any $N_1,N_2\in[1,\iy)$ with $N_1\le N_2$ we
have that
$$\sum_{p=N_1}^{N_2}\f{|\fr_{\bap}[\rM,\bap](p)|^2}{p^b}\ll\begin{cases}
N_1^{a-b}+N_2^{a-b}, & a\neq b,\\
\log N_1+\log N_2, & a=b.
\end{cases}$$ In particular, if $b=0$, then
$$\sum_{p=N_1}^{N_2}|\fr_{\bap}[\rM,\bap](p)|^2=\rR[\rM,\bap](N_2)-\rR[\rM,\bap](N_1-1).$$
\end{lemma}

\begin{thm} Let $\bap\in\BR^n$
be a vector of diophantine type $\kappa<\f{n-1}{n-2}$ $($$n\ge 2$$)$
such that the components of the vector $(\bap,1)$ are linearly
independent over $\BQ$, and let $\fA[\widehat\rM,\bap]$ be as
defined in $(1.7)$ where $\rM\in\fS^+_n(\BZ)$. If $K\le T^H$ for
some $H>0$, then we have that $\la|F^{K,\el}_{\rM}|^2\ra_T\ll
T^{-2\el}$ for any $\el\in\BN$ and moreover
\begin{equation}\lim_{K\to\iy}\la|F^{K,0}_{\rM}|^2\ra_T=\f{1}{2\pi^2}\,\fA[\widehat\rM,\bap].
\end{equation}
\end{thm}

\pf Let $\mu_{\el}(t)=\mu(t)\,t^{-2\el}$. Then it follows from
(2.10) that for each $\el\in\BN$,
\begin{equation}\begin{split}\la|F^{K,\el}_{\rM}|^2\ra_T&=\f{P^2_{\el}(\det\rM)^{\f{n-1+2\el}{2}}}{T^{2\el}}
\sum_{p=1}^{K^{2+\zt}}\sum_{q=1}^{K^{2+\zt}}\f{\fr[\widehat\rM,\bap](p)\,\overline{\fr[\widehat\rM,\bap](q)}}{(pq)^{\f{n+1+2\el}{4}}}\\
&\qquad\qquad\times\widehat \vp\biggl(\f{\sqrt
p}{K(\det\rM)^{\f{1}{2}}}\biggr)\overline{\widehat
\vp\biggl(\f{\sqrt
q}{K(\det\rM)^{\f{1}{2}}}\biggr)}\\&\quad\times\f{1}{T}\int_0^{\iy}\,\cos\biggl(\displaystyle\f{2\pi
t\sqrt
p}{(\det\rM)^{\f{1}{2}}}+\phi_{\el}\biggr)\cos\biggl(\displaystyle\f{2\pi
t\sqrt
q}{(\det\rM)^{\f{1}{2}}}+\phi_{\el}\biggr)\mu\biggl(\f{t}{T}\biggr)\f{T^{2\el}}{t^{2\el}}\,dt\\
&=\f{P^2_{\el}(\det\rM)^{\f{n-1-2\el}{2}}}{T^{2\el}}
\sum_{p=1}^{K^{2+\zt}}\sum_{q=1}^{K^{2+\zt}}\f{\fr[\widehat\rM,\bap](p)\,\overline{\fr[\widehat\rM,\bap](q)}}{(pq)^{\f{n+1+2\el}{4}}}\\
&\qquad\qquad\times\widehat \vp\biggl(\f{\sqrt
p}{K(\det\rM)^{\f{1}{2}}}\biggr)\overline{\widehat
\vp\biggl(\f{\sqrt q}{K(\det\rM)^{\f{1}{2}}}\biggr)}\\
&\quad\times\biggl(\f{1}{2}\,\RE\biggl[e^{-2i\phi_{\el}}\widehat\mu_{\el}\biggl(\f{T(\sqrt
p+\sqrt
q)}{(\det\rM)^{\f{1}{2}}}\biggr)\biggr]+\f{1}{2}\,\RE\biggl[\widehat\mu_{\el}\biggl(\f{T(\sqrt
p-\sqrt q)}{(\det\rM)^{\f{1}{2}}}\biggr)\biggr]\biggr)\\
&:=G^{K,T,\el}_{\rM} +H^{K,T,\el}_{\rM},
\end{split}\end{equation}
where $G^{K,T,\el}_{\rM}$ and $H^{K,T,\el}_{\rM}$ are the sums of
terms containing $\mu_{\el}\biggl(\f{T(\sqrt p+\sqrt
q)}{(\det\rM)^{\f{1}{2}}}\biggr)$ and $\mu_{\el}\biggl(\f{T(\sqrt
p-\sqrt q)}{(\det\rM)^{\f{1}{2}}}\biggr)$, respectively. Since
$|\widehat\mu_{\el}(s)|\ll(1+|s|)^{-N}$ for any $N>0$, and
$$\biggl|\widehat\mu_{\el}\biggl(\f{T(\sqrt
p+\sqrt
q)}{(\det\rM)^{\f{1}{2}}}\biggr)\biggr|\ll(\det\rM)^{\f{N}{2}}T^{-N}p^{-\f{N}{4}}q^{-\f{N}{4}},$$
it follows from Property 1 and Lemma 4.3 that the first term
$G^{K,T,\el}_{\rM}$ is bounded by
\begin{equation}|P_{\el}|^2(\det\rM)^{\f{n-1-2\el+N}{2}}T^{-2\el-N}
\sum_{p=1}^{K^{2+\zt}}\f{|\fr[\widehat\rM,\bap](p)|^2}{p^{\f{n+1+2\el}{2}+\f{N}{4}}}
\sum_{q=1}^{K^{2+\zt}}\f{1}{q^{\f{N}{4}}}\ll
T^{-2\el-N}.\end{equation} This implies that $G^{K,T,\el}_{\rM}$ has
a negligible contribution to $\la|F^{K,\el}_{\rM}|^2\ra_T$.

Next we estimate the part $H^{K,T,\el}_{\rM}$. We now split this sum
into the diagonal part $D^{K,T,\el}_{\rM}$ and the off-diagonal part
$O^{K,T,\el}_{\rM}$, that is,
$H^{K,T,\el}_{\rM}=D^{K,T,\el}_{\rM}+O^{K,T,\el}_{\rM}$ where
\begin{equation}D^{K,T,\el}_{\rM}=\f{\widehat\mu_{\el}(0)}{2}\f{P^2_{\el}(\det\rM)^{\f{n-1-2\el}{2}}}{T^{2\el}}
\sum_{p=1}^{K^{2+\zt}}\f{|\fr[\widehat\rM,\bap](p)|^2}{p^{\f{n+1+2\el}{2}}}\,
\biggl|\,\widehat\vp\biggl(\f{\sqrt
p}{K(\det\rM)^{\f{1}{2}}}\biggr)\biggr|^2
\end{equation}
and
\begin{equation}\begin{split}O^{K,T,\el}_{\rM}&=\f{1}{2}\f{P^2_{\el}(\det\rM)^{\f{n-1-2\el}{2}}}{T^{2\el}}
\sum_{\substack{p,q=1\\p\neq q}}^{K^{2+\zt}}\f{\fr[\widehat\rM,\bap](p)\,\overline{\fr[\widehat\rM,\bap](q)}}{(pq)^{\f{n+1+2\el}{4}}}\\
&\qquad\times\widehat\vp\biggl(\f{\sqrt
p}{K(\det\rM)^{\f{1}{2}}}\biggr)\overline{\widehat
\vp\biggl(\f{\sqrt
q}{K(\det\rM)^{\f{1}{2}}}\biggr)}\RE\biggl[\widehat\mu_{\el}\biggl(\f{T(\sqrt
p-\sqrt q)}{(\det\rM)^{\f{1}{2}}}\biggr)\biggr]
\end{split}\end{equation}
for $\el=0,1,\cdots\,.$ Here we note that $D^{K,T,\el}_{\rM}$
depends only on $K$, and so we write
$D^{K,T,0}_{\rM}=D^{K,0}_{\rM}$.

In order to prove Theorem 4.4, we show the following two lemmas on
the estimates for the diagonal part $D^{K,T,\el}_{\rM}$ and the
off-diagonal part $O^{K,T,\el}_{\rM}$.

\begin{lemma} Let $\bap\in\BR^n$ be a
vector of diophantine type $\kappa<\f{n-1}{n-2}$ $($$n\ge 2$$)$ such
that the components of the vector $(\bap,1)$ are linearly
independent over $\BQ$. If $\,\rM\in\fS^+_n(\BZ)$ and $K\le T^H$ for
some $H>0$, then we have that $D^{K,T,\el}_{\rM}\ll T^{-2\el}$ for
any $\el\in\BN$ and moreover
\begin{equation}\lim_{K\to\iy}D^{K,T,0}_{\rM}=\f{1}{2\pi^2}\,\fA[\widehat\rM,\bap].
\end{equation}
\end{lemma}

\pf The estimate $D^{K,T,\el}_{\rM}\ll T^{-2\el}$ for any
$\el\in\BN$ can be derived from the boundedness of the function
$\widehat g_{\psi}$ and Lemma 4.3. For the proof of (4.9), we recall
that $\widehat\mu(0)=1$, $P_0^2=1/\pi^2$ and $\widehat
g_{\psi}(0)=1$ by (2.1). Hence the required result follows from the
Lebesgue convergence theorem. \qed

\begin{lemma} Let $\bap\in\BR^n$ be a
vector of diophantine type $\kappa<\f{n-1}{n-2}$ $($$n\ge 2$$)$ such
that the components of the vector $(\bap,1)$ are linearly
independent over $\BQ$. If $\,\rM\in\fS^+_n(\BZ)$ and $K\le T^H$ for
some $H>0$, then we have that $\,\,O^{K,T,\el}_{\rM}\ll
T^{-2\el-\dt}$ for any $\el\ge 0$ and any $\dt\in(0,1)$.
\end{lemma}

\pf Take any $T>0$ large enough and any $\dt\in(0,1)$. We decompose
the sum of the off-diagonal part $O^{K,T,\el}_{\rM}$ into two sums
as follows; $O^{K,T,\el,1}_{\rM}$ for $|\sqrt p-\sqrt q|\ge
T^{-\dt}$ and $O^{K,T,\el,2}_{\rM}$ for $|\sqrt p-\sqrt
q|<T^{-\dt}$. As in Lemma 2.7 \cite{KS}, the estimates
$O^{K,T,\el,1}_{\rM}\ll T^{-2\el-\dt}$ and $O^{K,T,\el,2}_{\rM}\ll
T^{-2\el-\dt}$ follow from Lemma 4.3 and Property 1. \qed

Therefore we  can easily complete the proof of Theorem 4.4 by
applying (4.6), Lemma 4.5 and Lemma 4.6.\qed

{\bf Proof of Theorem 1.4.} It can be shown from Theorem 4.4, Lemma
4.5 and Lemma 4.6.\qed

\section{ Proof of Theorem 1.5}

Throughout this section, suppose that Property 1 holds.

We write again (1.2) as
\begin{equation}S_{\rM}(t,\e)=\f{F_{\rM}(t+\e)-F_{\rM}(t)}{\sqrt\e}+P_{\rM}(t,\e),
\end{equation}
where $P_{\rM}(t,\e):=\ds\f{1}{\sqrt\e}
\biggl(\f{(t+\e)^{\f{n-1}{2}}}{t^{\f{n-1}{2}}}-1\biggr)F_{\rM}(t+\e)$.
We mention an elementary observation used in \cite{KS}: there are
some $t_0,t_1\in (0,\iy)$ with $t_0<t_1$ such that $\mu\equiv 0$
outside the interval $[t_0,t_1]$. Then we define a non-negative
function $\upsilon\in C^{\iy}_c(\BR)$ such that
\begin{equation*}\upsilon(t)=\begin{cases}\sup_{[t_0,t_1]}\mu, & t\in[t_0/2,2 t_1],\\
0, & t\le t_0/4\text{ or }t\ge 4 t_1.\end{cases}\end{equation*} If
$\e$ is sufficiently small, then for any function $f\in L^1_{\rm
loc}(\BR)$ we have that
\begin{equation}\la|f(\,\cdot+\e)|\ra_{\mu,T}\le\|\upsilon\|_{L^1(\BR)}\la|f|\ra_{\om,T},
\end{equation}
where $\om=\upsilon/\|\upsilon\|_{L^1(\BR)}$. The estimate of
$S_{\rM}(\,\cdot\,,\e)$ heavily relies on the asymptotic value of
the normalized deviation $S_{\rM}(t,\e;K)$ of $F^{K,0}_{\rM}(t)$
given by
$$S_{\rM}(t,\e;K):=\f{F^{K,0}_{\rM}(t+\e)-F^{K,0}_{\rM}(t)}{\sqrt\e}.$$

\begin{lemma} Let $\bap\in\BR^n$ be a
vector of diophantine type $\kappa<\f{n-1}{n-2}$ $($$n\ge 2$$)$ such
that the components of the vector $(\bap,1)$ are linearly
independent over $\BQ$. If $\,\rM\in\fS^+_n(\BZ)$ and $K=T^n$ and
$T^{-\sm}\ll\e$ for some $\sm\in(0,1)$, then we have that
$$\lim_{\e\to
0}\la|S_{\rM}(\,\cdot\,,\e)-S_{\rM}(\,\cdot\,,\e;K)|^2\ra_T=0.$$
\end{lemma}

\pf  By Taylor's expansion, we have that
$$P_{\rM}(t,\e)=\f{\sqrt\e}{t}\,F_{\rM}(t+\e)\,\kappa(\e,t,n)$$ where
$\kappa(\e,t,n)=\sum_{k=1}^{\iy}\f{1}{k!}\f{n-1}{2}(\f{n-1}{2}-1)\cdots
(\f{n-1}{2}-k+1)\,\e^k t^{-k}$. Here we note that
$\kappa(\e,t,n)=\cO(1)$ for any sufficiently small $\e>0$. Since
$\mu_T$ is supported in $[t_0 T,t_1 T]$ (i.e. $0<t_0 T\le t\le t_1
T<\iy$), it follows from (5.2) and Theorem 1.4 that
\begin{equation*}\la|P_{\rM}(\,\cdot\,,\e)|^2\ra_T\ll\f{\e}{T^2}\,\|\upsilon\|_{L^1(\BR)}\la|F_{\rM}|^2\ra_{\om,T}
\to 0\,\,\,\,\,\text{as $\e\to 0$.}
\end{equation*}
Thus by (5.1) we obtain that
\begin{equation*}\lim_{\e\to
0}\bigl\la\bigl|S_{\rM}(\,\cdot\,,\e)-\f{F_{\rM}(\,\cdot+\e)-F_{\rM}(t)}{\sqrt\e}\bigr|^2\bigr\ra_T=0.
\end{equation*}
Since $K=T^n$, by Lemma 4.2 we see that
\begin{equation*}\f{1}{\e}\,\la|F_{\rM}-F^K_{\rM}|^2\ra_T\ll\f{1}{\e\,T}\ll\f{1}{T^{1-\sm}}\to
0\,\,\,\,\text{ as $T\to\iy$.}
\end{equation*}
Hence by (5.2) we have that
\begin{equation*}\f{1}{\e}\la|F_{\rM}(\,\cdot+\e)-F^K_{\rM}(\,\cdot+\e)|^2\ra_T\le\f{1}{\e}\,
\|\upsilon\|_{L^1(\BR)}\la|F_{\rM}-F^K_{\rM}|^2\ra_{\om,T}\to
0\,\,\,\,\,\text{as $\e\to 0$.}
\end{equation*}
This implies that
\begin{equation*}\lim_{\e\to
0}\bigl\la\bigl|S_{\rM}(\,\cdot\,,\e)-\f{F^K_{\rM}(\,\cdot+\e)-F^K_{\rM}(t)}{\sqrt\e}\bigr|^2\bigr\ra_T=0.
\end{equation*}
By applying (2.8), (2.9) and Theorem 4.4, we now proceed the same
type of estimates as above to obtain that
\begin{equation*}\lim_{\e\to
0}\bigl\la\bigl|S_{\rM}(\,\cdot\,,\e)-\f{F^{K,0}_{\rM}(\,\cdot+\e)-F^{K,0}_{\rM}(t)}{\sqrt\e}\bigr|^2\bigr\ra_T=0.
\end{equation*} Therefore we complete the proof. \qed

We now give the explicit formula of
$\la|S_{\rM}(\,\cdot\,,\e;K)|^2\ra_T$ by the help of the previous
arguments. From (2.10), we have that
\begin{equation*}\begin{split}&F^{K,0}_{\rM}(t+\e)-F^{K,0}_{\rM}(t)=\f{(\det\rM)^{\f{n-1}{4}}}{\pi}\\
&\qquad\qquad\times\sum_{p=1}^{K^{2+\zt}}\,\biggl[\cos\biggl(\displaystyle\f{2\pi
(t+\e)\sqrt
p}{(\det\rM)^{\f{1}{2}}}+\phi_0\biggr)-\cos\biggl(\displaystyle\f{2\pi
t\sqrt
p}{(\det\rM)^{\f{1}{2}}}+\phi_0\biggr)\biggr]\\
&\qquad\qquad\quad\qquad\qquad\times\f{\fr[\widehat\rM,\bap](p)}{p^{\f{n+1}{4}}}\,\widehat
\vp\biggl(\f{\sqrt p}{K(\det\rM)^{\f{1}{2}}}\biggr)\\
&\qquad\qquad=\f{-2(\det\rM)^{\f{n-1}{4}}}{\pi}\sum_{p=1}^{K^{2+\zt}}\,\sin\biggl(\displaystyle\f{2\pi
(t+\e/2)\sqrt
p}{(\det\rM)^{\f{1}{2}}}+\phi_0\biggr)\sin\biggl(\displaystyle\f{\pi\e\sqrt
p}{(\det\rM)^{\f{1}{2}}}\biggr)\\
&\qquad\qquad\quad\qquad\qquad\times\f{\fr[\widehat\rM,\bap](p)}{p^{\f{n+1}{4}}}\,\widehat
\vp\biggl(\f{\sqrt p}{K(\det\rM)^{\f{1}{2}}}\biggr)
\end{split}\end{equation*}
where $\phi_0=-\f{n+1}{4}\pi$. Thus this yields that
\begin{equation}\begin{split}&\la|S_{\rM}(\,\cdot\,,\e;K)|^2\ra_T=\f{4(\det\rM)^{\f{n-1}{2}}}{\e\pi^2}\sum_{p,q=1}^{K^{2+\zt}}
\sin\biggl(\displaystyle\f{\pi\e\sqrt
p}{(\det\rM)^{\f{1}{2}}}\biggr)\sin\biggl(\displaystyle\f{\pi
\e\sqrt q}{(\det\rM)^{\f{1}{2}}}\biggr)\\
&\qquad\qquad\quad\times\f{\fr[\widehat\rM,\bap](p)\,\overline{\fr[\widehat\rM,\bap](q)}}{(pq)^{\f{n+1}{4}}}\,\widehat
\vp\biggl(\f{\sqrt
p}{K(\det\rM)^{\f{1}{2}}}\biggr)\overline{\widehat
\vp\biggl(\f{\sqrt
q}{K(\det\rM)^{\f{1}{2}}}\biggr)}\\
&\qquad\qquad\quad\qquad\times\bigl\la\sin\bigl(\theta_{\rM}(\,\cdot\,,\e,\sqrt
p)+\phi_0\bigl) \sin\bigl(\theta_{\rM}(\,\cdot\,,\e,\sqrt
q)+\phi_0\bigr)\bigr\ra_T,
\end{split}\end{equation}
where
$$\theta_\rM(t,\e,s)=\displaystyle\f{2\pi(t+\e/2)s}{(\det\rM)^{\f{1}{2}}}.$$
We note that
\begin{equation}\begin{split}&\bigl\la\sin\bigl(\theta_{\rM}(\,\cdot\,,\e,\sqrt
p)+\phi_0\bigl) \sin\bigl(\theta_{\rM}(\,\cdot\,,\e,\sqrt
q)+\phi_0\bigr)\bigr\ra_T\\
&\qquad\qquad\qquad=-\f{1}{2}\,\RE\biggl[e^{i\theta_{\rM}(0,\e,\sqrt
p+\sqrt q)+i 2\phi_0}\,\widehat\mu\biggl(\f{-T(\sqrt p+\sqrt
q)}{(\det\rM)^{1/2}}\biggr)\biggr]\\
&\qquad\qquad\qquad\qquad+\f{1}{2}\,\RE\biggl[e^{-i\theta_{\rM}(0,\e,\sqrt
p-\sqrt q)}\,\widehat\mu\biggl(\f{T(\sqrt p-\sqrt
q)}{(\det\rM)^{1/2}}\biggr)\biggr].
\end{split}\end{equation}
Since $|\widehat\mu|$ has fast decay, as in the proof of Theorem 4.4
the first quantity in (5.4) has a negligible contribution to
$\la|S_{\rM}(\,\cdot\,,\e;K)|^2\ra_T$ as $T\to\iy$. Now the
remaining sum of $\la|S_{\rM}(\,\cdot\,,\e;K)|^2\ra_T$ can be
decomposed into the diagonal part $S^D_{\rM}(\e,K)$ and the
off-diagonal part $S^O_{\rM}(\e,K,T)$ as follows;
\begin{equation}\begin{split}S^D_{\rM}(\e,K)&=\f{2(\det\rM)^{\f{n-1}{2}}}{\e\pi^2}\sum_{p=1}^{K^{2+\zt}}
\sin^2\biggl(\displaystyle\f{\pi\e\sqrt
p}{(\det\rM)^{\f{1}{2}}}\biggr)\\
&\qquad\qquad\qquad\qquad\qquad\times\f{|\fr[\widehat\rM,\bap](p)|^2}{p^{\f{n+1}{2}}}\,\biggl|\,\widehat
\vp\biggl(\f{\sqrt p}{K(\det\rM)^{\f{1}{2}}}\biggr)\biggr|^2
\end{split}\end{equation} and
\begin{equation}\begin{split}S^O_{\rM}(\e,K,T)&=\f{2(\det\rM)^{\f{n-1}{2}}}{\e\pi^2}\sum_{\substack{p,q=1\\p\neq q}}^{K^{2+\zt}}
\sin\biggl(\displaystyle\f{\pi\e\sqrt
p}{(\det\rM)^{\f{1}{2}}}\biggr)\sin\biggl(\displaystyle\f{\pi
\e\sqrt q}{(\det\rM)^{\f{1}{2}}}\biggr)\\
&\times\f{\fr[\widehat\rM,\bap](p)\,\overline{\fr[\widehat\rM,\bap](q)}}{(pq)^{\f{n+1}{4}}}\,\widehat
\vp\biggl(\f{\sqrt
p}{K(\det\rM)^{\f{1}{2}}}\biggr)\overline{\widehat
\vp\biggl(\f{\sqrt q}{K(\det\rM)^{\f{1}{2}}}\biggr)}\\
&\times\RE\biggl[e^{-i\theta_{\rM}(0,\e,\sqrt p-\sqrt
q)}\,\widehat\mu\biggl(\f{T(\sqrt p-\sqrt
q)}{(\det\rM)^{1/2}}\biggr)\biggr].
\end{split}\end{equation}
Here we note that the diagonal part $S^D_{\rM}(\e,K)$ is independent
of $T$.

As in the proof of Theorem 1.4, we show that the off-diagonal part
$S^O_{\rM}(\e,K,T)$ makes no contribution to
$\la|S_{\rM}(\,\cdot\,,\e;K)|^2\ra_T$ and we obtain the explicit
asymptotics of the diagonal part of $S^D_{\rM}(\e,K)$.

\begin{lemma} Let $\bap\in\BR^n$ be a
vector of diophantine type $\kappa<\f{n-1}{n-2}$ $($$n\ge 2$$)$ such
that the components of the vector $(\bap,1)$ are linearly
independent over $\BQ$. If $\,\rM\in\fS^+_n(\BZ)$ and $K\le T^H$ for
some $H>0$, then we have that $$S^O_{\rM}(\e,K,T)\ll
\e^{-1}T^{-\dt}$$ for any $\e>0$ and any $\dt\in(0,1)$.
\end{lemma}

\pf Take any $\e>0$. As in the proof of Lemma 4.6 with $\el=0$, we
then obtain the estimate $$\e\,S^O_{\rM}(\e,K,T)\ll T^{-\dt}$$ for
any $\dt\in(0,1)$.\qed

\begin{lemma} Let $\bap\in\BR^n$ be a
vector of diophantine type $\kappa<\f{n-1}{n-2}$ $($$n\ge 2$$)$ such
that the components of the vector $(\bap,1)$ are linearly
independent over $\BQ$. If $\,\rM\in\fS^+_n(\BZ)$, then for any
$N_1,N_2\in(0,\iy)$ with $N_1<N_2$, we have that
$$\lim_{\e\to 0}\sum_{N_1\le\e^2
p<N_2}\e^n\,|\fr[\rM,\bap](p)|^2=|E^{\rM}|(N_2^{\f{n}{2}}-N_1^{\f{n}{2}}).$$
\end{lemma}

\pf It easily follows from Property 1 and Lemma 4.3 that
\begin{equation}\begin{split}&\limsup_{\e\to 0}\sum_{N_1\le\e^2
p<N_2}\e^n\,|\fr[\rM,\bap](p)|^2\\&\qquad\qquad\le\lim_{\e\to
0}\e^n\bigl(\rR[\rM,\bap]([N_2\e^{-2}])-\rR[\rM,\bap]([N_1\e^{-2}])\bigr)\\
&\qquad\qquad=|E^{\rM}|(N_2^{\f{n}{2}}-N_1^{\f{n}{2}}).
\end{split}\end{equation}
Also the lower bound can be obtained in a similar way. \qed

Given some fixed numbers $\vep$ and $N$ with $0<\vep<N<\iy$, we
consider a truncated sum $\underline S^{D,0}_{\rM}(\vep,N,K;\e)$
related with $S^D_{\rM}(\e,K)$ defined by
\begin{equation}\underline S^{D,0}_{\rM}(\vep,N,K;\e)=\f{2(\det\rM)^{\f{n-1}{2}}}{\e\pi^2}\sum_{\vep\e^{-2}\le p<N\e^{-2}}
\sin^2\biggl(\displaystyle\f{\pi\e\sqrt
p}{(\det\rM)^{\f{1}{2}}}\biggr)\f{|\fr[\widehat\rM,\bap](p)|^2}{p^{\f{n+1}{2}}}.
\end{equation}

\begin{lemma} Let $\bap\in\BR^n$ be a
vector of diophantine type $\kappa<\f{n-1}{n-2}$ $($$n\ge 2$$)$ such
that the components of the vector $(\bap,1)$ are linearly
independent over $\BQ$. If $\,\rM\in\fS^+_n(\BZ)$, then for any
$\vep,N\in(0,\iy)$ with $\vep<N$, we have that
$$\lim_{\e\to
0}\underline
S^{D,0}_{\rM}(\vep,N,K;\e)=\f{n|E^{\rM}|}{\pi^2}\int_{\vep_0}^{N_0}
\f{\sin^2(\pi\sqrt t)}{t^{\f{3}{2}}}\,dt$$ where
$\vep_0=\vep/\det\rM$ and $N_0=N/\det\rM$.
\end{lemma}

\pf Split the interval $[\vep,N)$ into $L$ equal intervals with
length $l=(N-\vep)/L$ and decompose $\underline
S^{D,0}_{\rM}(\vep,N,K;\e)$ into $\underline
S^{D,0,k}_{\rM}(\vep,N,K;\e)$ ($k=0,1,\cdots,L-1$) defined by
$$\underline S^{D,0,k}_{\rM}(\vep,N,K;\e)=\f{2(\det\rM)^{\f{n-1}{2}}}{\e\pi^2}\sum_{\vep+lk\le\e^2 p<\vep+l(k+1)}
\sin^2\biggl(\displaystyle\f{\pi\e\sqrt
p}{(\det\rM)^{\f{1}{2}}}\biggr)\f{|\fr[\widehat\rM,\bap](p)|^2}{p^{\f{n+1}{2}}}.$$
For $t>0$, we set $F(t)=t^{-\f{n+1}{2}}\sin^2(\pi\sqrt t).$ Then we
have that
$$\underline S^{D,0,k}_{\rM}(\vep,N,K;\e)=\f{2\,\e^n}{\pi^2(\det\rM)}\sum_{\vep+lk\le\e^2 p<\vep+l(k+1)}
F\biggl(\f{\e^2 p}{\det\rM}\biggr)|\fr[\widehat\rM,\bap](p)|^2.$$ We
set
$R^{D,0}_{\rM}(\vep,N,K;\e,l)=\sum_{k=0}^{L-1}R^{D,0,k}_{\rM}(\vep,N,K;\e,l)$
where
$$R^{D,0,k}_{\rM}(\vep,N,K;\e,l)=\f{2\,\e^n}{\pi^2(\det\rM)}\,F\biggl(\f{\vep+lk}{\det\rM}\biggr)\sum_{\vep+lk\le\e^2
p<\vep+l(k+1)}|\fr[\widehat\rM,\bap](p)|^2.$$ Since $|F'(t)|\ll
t^{-\f{n+1}{2}}$ for $t>0$, by the mean value theorem we have that
$$\biggl|F\biggl(\f{\e^2 p}{\det\rM}\biggr)-F\biggl(\f{\vep+lk}{\det\rM}\biggr)\biggr|
\ll l(\vep+lk)^{-\f{n+1}{2}}(\det\rM)^{\f{n-1}{2}}$$ and
\begin{equation}\begin{split}(\vep+l(k+1))^{\f{n}{2}}-(\vep+lk)^{\f{n}{2}}&=\f{nl}{2}\int_0^1(\vep+l(k+s))^{\f{n}{2}-1}ds\\
&\le\f{nl}{2}(\vep+l(k+1))^{\f{n}{2}-1}.
\end{split}\end{equation}
From the fact that
$|E^{\rM}|=(\det\rM)^{\f{n}{2}-1}|E^{\widehat\rM}|$, it follows from
(5.7) that
\begin{equation}\begin{split}&|\underline S^{D,0}_{\rM}(\vep,N,K;\e)-R^{D,0}_{\rM}(\vep,N,K;\e,l)|\\
&\ll\f{2(\det\rM)^{\f{n-3}{2}}l}{\pi^2}\sum_{k=0}^{L-1}(\vep+lk)^{-\f{n+1}{2}}\sum_{\vep+lk\le\e^2
p<\vep+l(k+1)}\e^n\,|\fr[\widehat\rM,\bap](p)|^2\\
&\le\f{n|E^{\rM}|\,l^2}{\pi^2(\det\rM)^{\f{1}{2}}}\sum_{k=0}^{L-1}\biggl(\f{\vep+l(k+1)}{\vep+lk}\biggr)^{\f{n+1}{2}}
\f{1}{(\vep+l(k+1))^{\f{3}{2}}}\\
&\le\f{2^{\f{n+1}{2}}n|E^{\rM}|\,l^{\f{1}{2}}
N_0}{\pi^2(\det\rM)^{\f{1}{2}}}\to 0\,\,\,\,\text{ as $\e\to 0$ and
$l\to 0$,}
\end{split}\end{equation}
where $N_0=\sum_{k=0}^{\iy}(k+1)^{-\f{3}{2}}<\iy$. Hence this
implies that
$$\lim_{\e\to 0}\underline S^{D,0}_{\rM}(\vep,N,K;\e)=\lim_{l\to 0}\lim_{\e\to
0}R^{D,0}_{\rM}(\vep,N,K;\e,l).$$ Now we estimate
$R^{D,0}_{\rM}(\vep,N,K;\e,l)$. By Lemma 5.3 and (5.9) we have that
\begin{equation}\begin{split}&\lim_{\e\to 0}R^{D,0}_{\rM}(\vep,N,K;\e,l)\\
&=\f{2}{\pi^2(\det\rM)}\sum_{k=0}^{L-1}F\biggl(\f{\vep+lk}{\det\rM}\biggr)\,\lim_{\e\to
0}\sum_{\vep+lk\le\e^2
p<\vep+l(k+1)}\e^n\,|\fr[\widehat\rM,\bap](p)|^2\\
&=\f{n(\det\rM)^{\f{n}{2}-1}|E^{\widehat\rM}|}{\pi^2}\sum_{k=0}^{L-1}
F\biggl(\f{\vep+lk}{\det\rM}\biggr)\f{l}{\det\rM}\int_0^1\biggl(\f{\vep+l(k+s)}{\det\rM}\biggr)^{\f{n}{2}-1}\,ds.
\end{split}\end{equation}
Using the facts that
$|E^{\rM}|=(\det\rM)^{\f{n}{2}-1}|E^{\widehat\rM}|$ and
$$\lim_{l\to 0}\int_0^1\biggl(\f{\vep+l(k+s)}{\vep+lk}\biggr)^{\f{n}{2}-1}ds=1\,\,\text{ for all $k=0,1,\cdots,L-1$,} $$
by (5.11) we obtain that
\begin{equation*}\begin{split}\lim_{\e\to
0}\underline S^{D,0}_{\rM}(\vep,N,K;\e)&=\lim_{l\to 0}\lim_{\e\to
0}R^{D,0}_{\rM}(\vep,N,K;\e,l)\\
&=\lim_{l\to 0}\f{n|E^{\rM}|}{\pi^2}\sum_{k=0}^{L-1}
F\biggl(\f{\vep+lk}{\det\rM}\biggr)\f{l}{\det\rM}\biggl(\f{\vep+lk}{\det\rM}\biggr)^{\f{n}{2}-1}\\
&=\lim_{l\to 0}\sum_{k=0}^{L-1}
G\biggl(\f{\vep+lk}{\det\rM}\biggr)\f{l}{\det\rM}=\int_{\vep_0}^{N_0}
G(t)\,dt,
\end{split}\end{equation*}
where
$G(t)=n|E^{\rM}|\pi^{-2}F(t)t^{\f{n}{2}-1}=n|E^{\rM}|\pi^{-2}t^{-\f{3}{2}}\sin^2(\pi\sqrt
t)$ for $t>0$. Since $G$ is continuous on $[\vep_0,N_0]$, this
completes the proof. \qed

\begin{remark} We observe that the estimates in Lemma
5.4 no longer depend on $\vep$ and $N$.
\end{remark}

\begin{lemma} Let $\bap\in\BR^n$ be a
vector of diophantine type $\kappa<\f{n-1}{n-2}$ $($$n\ge 2$$)$ such
that the components of the vector $(\bap,1)$ are linearly
independent over $\BQ$. If $\,\rM\in\fS^+_n(\BZ)$, then we have that
$\,\lim_{\e\to 0}S^D_{\rM}(\e,K)=n|E^{\rM}|\,$ as $\e\to 0$ and $\e
K\to\iy$.
\end{lemma}

\pf Take any $\vep,N\in(0,\iy)$ with $\vep<N$. We split
$S^D_{\rM}(\e,K)$ given in (5.5) into three parts as follows;
\begin{equation*}\begin{split}S^D_{\rM}(\e,K)&=S^{D-}_{\rM}(\vep,K;\e)+S^{D0}_{\rM}(\vep,N,K;\e)+S^{D+}_{\rM}(N,K;\e)\\
&:=\f{2(\det\rM)^{\f{n-1}{2}}}{\e\pi^2}\sum_{1\le p<\vep\e^{-2}}
H_{\bap}(\e,K,\rM,p,\psi)\\
&\qquad+\f{2(\det\rM)^{\f{n-1}{2}}}{\e\pi^2}\sum_{\vep\e^{-2}\le
p<N\e^{-2}}H_{\bap}(\e,K,\rM,p,\psi)\\
&\qquad+\f{2(\det\rM)^{\f{n-1}{2}}}{\e\pi^2}\sum_{N\e^{-2}\le p\le
K^{2+\zt}} H_{\bap}(\e,K,\rM,p,\psi),
\end{split}\end{equation*} where $H_{\bap}(\e,K,\rM,p,\psi)$ is
given by
$$H_{\bap}(\e,K,\rM,p,\psi)=\sin^2\biggl(\displaystyle\f{\pi\e\sqrt
p}{(\det\rM)^{\f{1}{2}}}\biggr)\f{|\fr[\widehat\rM,\bap](p)|^2}{p^{\f{n+1}{2}}}
\biggl|\,\widehat\vp\biggl(\f{\sqrt
p}{K(\det\rM)^{\f{1}{2}}}\biggr)\biggr|^2.$$ Using the fact that
$\sin^2(t)\le t^2$ for small $t>0$ and $|\sin^2(t)|\le 1$, it
follows from Property 1 and Lemma 4.3 that
\begin{equation}\begin{split}S^{D-}_{\rM}(\vep,K;\e)&\ll\f{(\det\rM)^{\f{n-1}{2}}}{\e\pi^2}\f{\pi^2\e^2
p}{\det\rM} \sum_{1\le
p<\vep\e^{-2}}\f{|\fr[\widehat\rM,\bap](p)|^2}{p^{\f{n+1}{2}}}\\&\ll
(\det\rM)^{\f{n-3}{2}}\sqrt\vep\end{split}\end{equation} and
\begin{equation}\begin{split}S^{D+}_{\rM}(N,M;\e)&\ll\f{(\det\rM)^{\f{n-1}{2}}}{\e\pi^2}\sum_{N\e^{-2}\le
p\le
M^{2+\zt}}\f{|\fr[\widehat\rM,\bap](p)|^2}{p^{\f{n+1}{2}}}\\&\ll\f{(\det\rM)^{\f{n-1}{2}}}{\pi^2\sqrt
N}.\end{split}\end{equation} Finally we compare
$S^{D0}_{\rM}(\vep,N,K;\e)$ with $\underline
S^{D0}_{\rM}(\vep,N,K;\e)$. Since $|\widehat
\vp(t)|^2\ll\f{1}{1+t^2}$, we have that
$$\biggl|\,\bigl|\,\widehat\vp\biggl(\f{\sqrt
p}{K(\det\rM)^{\f{1}{2}}}\biggr)\bigr|^2-1\biggr|\ll\f{N}{\e^2
K^2(\det\rM)}.$$ Thus by applying Property 1 and Lemma 4.3 again we
obtain that
\begin{equation}\begin{split}&|S^{D0}_{\rM}(\vep,N,K;\e)-\underline S^{D0}_{\rM}(\vep,N,K;\e)|\le\f{2N(\det\rM)^{\f{n-3}{2}}}{\pi^2(\e
K)^2}\biggl(\f{1}{\sqrt\vep}+\f{1}{\sqrt N}\biggr)\to 0
\end{split}\end{equation}
as $\e K\to\iy$. Therefore we conclude that
$$\lim_{\e\to 0}S^D_{\rM}(\e,K)=\f{n|E^{\rM}|}{\pi^2}\int_0^{\iy}\f{\sin^2(\pi\sqrt t)}{t^{\f{3}{2}}}\,dt=n|E^{\rM}|$$
by combining (5.12), (5.13), (5.14) and Lemma 5.4. \qed

{\bf Proof of Theorem 1.5.} Let $K=T^n$ and let $\e\gg T^{-\gm}$ for
some $\gm\in(0,1)$ Then we see that $\e K\to\iy$ and $T\to\iy$ as
$\e\to 0$. From Lemma 5.6 and Lemma 5.2, we can obtain that
$\lim_{T\to\iy}\la|S_{\rM}(\,\cdot\,,\e;K)|\ra_T=n|E^{\rM}|.$
Therefore we complete the proof by applying Lemma 5.1.\qed

\noindent{\bf Acknowledgement.} We would like to thank Alexander V.
Sobolev for several helpful comments.

\end{document}